%% file: main.tex
\documentclass[preprint,12pt]{elsarticle}
\usepackage{amssymb}
\usepackage{amsthm}
\usepackage{lineno}

\usepackage{textcmds}
\usepackage[colorlinks]{hyperref}
\usepackage{caption}
\usepackage{tikz,siunitx}
\usetikzlibrary{shapes.geometric,shapes.symbols}
\newcommand{\tikzsymbol}[2][circle]{\tikz[baseline=-0.5ex]\node[inner
sep=2pt,shape=#1,draw,#2]{};}

\usepackage{comment}
\usepackage{lmodern}
\usepackage{graphicx}
\usepackage{amsmath}
\usepackage{amsthm}
\usepackage{mathtools}
\usepackage{cleveref}
\usepackage[caption=false]{subfig}
\usepackage{pgfplots}
\pgfplotsset{compat=1.18}

\usepackage{booktabs}
\theoremstyle{definition}

\theoremstyle{remark}

\usepackage{color}
\usepackage{sectsty}
\usepackage[percent]{overpic}

\graphicspath{ {./figures/} }

\begin{document}
\input{frontmatter.tex}


\input{mainbody.tex}

\bibliographystyle{elsarticle-num-names}
\bibliography{valve}

\newpage

\appendix

\input{appendix.tex}

\end{document}

%% file: frontmatter.tex
\begin{frontmatter}



\title{Enhanced stability of pressure relief valves: mechanistic design and analysis}


\author[inst1]{Hong Tang}

\affiliation[inst1]{organization={Department of Engineering Mathematics,University of Bristol},
            addressline={Queen's Building},
            city={Bristol},
            postcode={BS1 8TH},
            country={UK}}

\author[inst2]{Istvan Erdodi}

\author[inst1]{Alan R. Champneys}


\author[inst2]{Csaba J. H{ő}s}

\affiliation[inst2]{organization=
	{Department of Hydrodynamic Systems, Budapest University of Technology and Economics},
            addressline={M{ű}egyetem rkp. 3.},
            city={Budapest},
            postcode={1111},
            country={Hungary}
}

\input{Abstract.tex}


\begin{graphicalabstract}
\vspace{-2cm}
\includegraphics[width=1.2\textwidth]{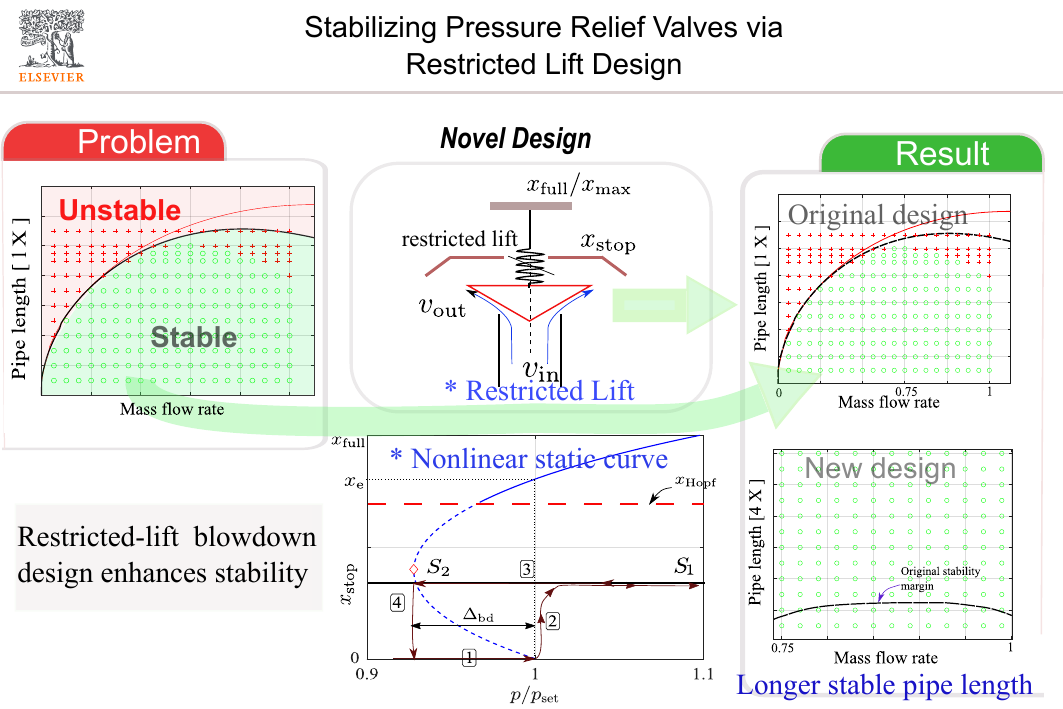}

\end{graphicalabstract}


\begin{highlights}
	\item A new design concept greatly enhances the stability of pressure-relief valves against flutter and chatter instabilities.
    \item The design exploits interaction between fold and Hopf bifurcations, and nonsmooth analysis of contacting pseudo equilibria. 
	\item Practical implementation is proposed either using oversized valves with restricted lift, or via novel equilibrium-versus-lift characteristics. 
	\item An experimentally validated coupled pipeline-valve model is improved to deal with the new design.
	 \item Detailed steady and transient simulations verify wide parameter sets for which oscillation-free operation can be achieved.
\end{highlights}

\begin{keyword}
process safety \sep nonlinear dynamics \sep fluid-structure  interaction \sep pressure-relief valve  \sep valve chatter \sep 
quarter-wave instability 
\PACS 47.85.L- \sep 46.32.+x
\MSC 35L35 \sep 37N10
\end{keyword}

\end{frontmatter}

%% file: Abstract.tex
\begin{abstract}
Pressure relief valves, often the critical last line of defense in process engineering, are known to be susceptible to valve chatter. Such behavior has been shown to arise from a 
flutter instability, or Hopf bifurcation, associated with the fundamental, quarter-wave acoustic mode of their inlet piping.  Here, a novel design concept is proposed and analyzed for eliminating this instability.
The concept involves using an oversized valve with reduced lift and adopting a discharge characteristic that enhances the blowdown effect, so that the valve opens immediately to its upper lift limit upon reaching set pressure.  
The concept is demonstrated numerically using an updated version of a 1D fluid pipe dynamics mathematical model that has been validated against experimental data. 
Stability properties are 
analysed using dynamical systems theory, applied to an improved reduced-order modal model.  
It is shown how the valve settles to a stable so-called pseudo equilibrium, in contact with the upper stop, provided the coefficient of restitution is not too large. Such stable operation is reached despite the equivalent regular valve being unstable to the quarter-wave Hopf bifurcation. Parameter studies using the reduced-order model demonstrate the extent of the enhanced stability effect, which is confirmed using the full  model for the case of gas service valves. 
\end{abstract}

%% file: mainbody.tex
%
%
\section{Introduction}\label{sec:introduction}

Pressure-relief valves (PRVs) are typically the last line of defense in pressurized equipment to protect vessels, pipelines, and pumps from exceeding safe operating limits. 
In the process industries, several major disasters over the last fifty years have been
attributed, at least in part, to improper maintenance or design of such valves and their inlet/outlet systems.
One of the major causes of such failure is a flutter instability, leading to damaging chattering behavior where the valve impacts with its seat at high frequency. Despite over 60 years of research, the prediction and alleviation of valve chatter remain unsolved problems. 

Specifically, we focus here only on the standard topology of direct spring operated pressure-relief valves (DSOPRV) commonly used in the oil and gas industry.  However the general principles we elucidate should be applicable to other kinds of relief valves  \cite{Cana_review} and more broadly to mechanical systems that experience instability through 
fluid-structure interaction.

\subsection{Chatter instability}

Since the origin of the spring-operated pressure-relief with the advent of steam engines, it has been widely understood that under certain conditions such valves can chatter. 
The first serious attempt to understand the causes
of chatter was the work of  Funk \cite{Funk} in the 1960s, 
who found that such valves are inclined to become unstable at a critical frequency that coincides
with the fundamental vibratory mode of the pipeline. Through subsequent experimental and theoretical studies by many authors over the following 50 years, it has  
been well-established that pressure-relief valve instability should not be considered in isolation, but that one needs to consider the dynamics of the whole system including the tank, valve and the inlet line connecting them. See, for example, the extensive list of references in  \cite{Darby:1,HosCh2014-I,Schmidt_models}.

A new urgency into finding a cure for chatter instabilities 
emerged in the 2010s through a series of scientific studies and industrial guidance. Most notable were the combined
numerical and experimental studies by Darby and co-workers
\cite{Darby:1,Darby:2,Darby:3}, detailed numerical work by  Melhem \cite{Melhem1,Melhem2} and the series of papers 
\cite{HosCh2014-I,HosCh2015-II,HosCh2016-III,HosCh2017-IV}.
In particular, in this latter series 
the authors developed a comprehensive theory of the so-called {\em quarter-wave} instability, backed up with experimental
measurements and numerical simulations (see also \cite{Ma_experiments,Schmidt_models} for more recent validation). 
It was found that all systems that involve a DSOPRV attached to inlet piping suffer from a fundamental flutter instability (or Hopf bifurcation to give it its mathematical name) in which the valve effectively applies negative damping to the fundamental quarter-wave
acoustic mode of the fluid in the inlet piping. 
The resulting oscillations rapidly grow into {\em chattering}, in which the valve repeatedly impacts with its seat at high frequency, causing damage or fatigue to the valve.
Ultimately, such behavior can prevent effective pressure relief.

The case of liquid service valves is more extreme than gas or vapor because the incompressibility of the fluid makes it liable to 
water-hammer effects
\cite{2010s:1,2010s:9,2010s:11} that can lead to downstream damage. Furthermore, in the liquid case, the underlying Hopf bifurcation is {\em subcritical} \cite{HosCh2016-III}, meaning that unstable valves transition directly to chatter without incipient flutter. 

Other attempts to understand the onset and extent of valve 
chatter have included detailed experimental studies \cite{2010s:6,2010s:8,2010s:12,Ma_experiments,new4,Smith2025},
full-scale computational fluid dynamics (CFD) of
the flow inside the valve \cite{Song2014,Wu2015,2010s:7,Erdodi_1}, simplified 1D-flow models \cite{Hos2014,2020:1,Safrany}, mechanical models of the valve-body impact itself \cite{2020:2,2020:4}, methods from nonlinear dynamical systems theory \cite{Izuchi,HosCh12gr,2010s:3,BazsoCh2014,BazsoHos2015,2010s:10}, delay loops
\cite{new3,Stepan:1,Stepan:2}
and various forms of analytic approximation \cite{2010s:2,Song2013,2010s:5,2010s:9,new7}.

Attempts have also been made to analyse instability mechanisms in more complex pressure-relief systems such as pilot-operated and other multi-stage valves 
\cite{2010s:4,2020:3,new8}.
The effects of valve chatter in specific applications have also been considered
including steam operations \cite{Domnick2018}, 
air suspension systems \cite{new2}, hydraulics \cite{new5}, 
subsea applications \cite{new6} and 
nuclear engineering \cite{new1}.

The conclusion seems to remain that the problem of chatter prediction is not solved. For example, 
Keszthelyi {\em et al.} \cite{Schmidt_models} after reviewing different modeling approaches conclude that \emph{`chattering of
safety valves is a rather stochastic process'} and that \emph{`only a more sophisticated approach
including difficult-to-measure parameters may allow for a more
reliable prediction of the stability of a safety valve'}.

Practical methods for avoiding pressure-relief-valve chatter also remain elusive.
Industry standards, such as API 520 part II \cite{API:520}, 
and ISO 4126 \cite{ISO:4126}
emphasize that to maintain reliable performance, the inlet line should usually be kept as short as 
possible so that pressure loss does not exceed 3\% of the valve’s set opening pressure. This criterion was originally proposed in the API code to avoid \emph{rapid cycling}; however, such problems are typically alleviated by designing valves with suitable blowdown (see Sec.~\ref{sec:blowdown} below). 

Nevertheless, as the API 520 points out, a shorter length $L$ of inlet piping is also advantageous for avoiding chatter. Specifically, direct calculation \cite{Hos2014,HosCh2015-II} shows that the quarter-wave instability will be triggered whenever the (dimensionless) mass flow rate $q$ through the valve is lower than a critical value, which increases with pipe length
$q<q_c(L)$, where to first approximation $ q_c(L) \propto L^{3/2}$. 

One solution would be to 
eliminate the inlet piping altogether and place PRVs directly on top of the vessels they are protecting.  This solution is not always practical and would
require significant redesign of existing valves where relief systems are often designed to a significant distance from
the pressure vessels so that they can discharge into a `safe' area. 
Also, even in the absence of inlet piping, DSOPRVs can still undergo chattering instabilities \cite{HosCh12gr}, which can be exacerbated when connected to outlet piping or bellows \cite{Chabane2009,Schmidt_outlet,Stepan:1,new3}.  

An alternative approach would be to 
greatly increase mechanical damping in the valve itself. At present, this is explicitly precluded in design codes such as API520 which state that the valve must \emph{`open unimpeded}'. Mechanical damping, such as frictional or viscous dampers, would increase the likelihood of the valve seizing after extended deployment in the field. In particular, in another recent paper 
\cite{Schmidt_2025}  Keszthelyi {\em et al.} experimentally tested a number of such possible damping solutions and found them to be impractical, ineffective or overly expensive, concluding
that \emph{`further research is necessary on how
such [damping] forces could be provided in a cost-effective manner.'}

Yet another idea is to employ damping in the fluid itself by allowing the inlet piping to be extremely long, so as to provide sufficient pipe friction to damp any acoustic modes. See e.g.~\cite{Izuchi} for some evidence of this working in theory. However, such a solution would be fraught with difficulty, 
because the amount of frictional damping depends on many unmeasurable flow and environmental factors, and in any case, such inlet piping would fail the 3\% rule. 

An alternative fundamental idea proposed to avoid chattering, and pertinent to the idea presented in this paper, is 
that valves should be designed to open and close rapidly \cite{Schmidt_models,HosCh2017-IV}, avoiding the low flow-rate regime where the quarter-wave  instability is most prevalent.  However for many valve and pipeline assemblies, instability can occur at the equilibrium lift position of the 
valve, which is itself typically a function of flow rate $q$. 
For a given flow rate and valve, the quarter-wave instability occurs upon increasing $L$.

\subsection{The new design concept}
\label{sec:design-concept-intro}

Here we shall briefly motivate and describe the key new idea in this paper. A detailed analysis will be presented in the main body of the text.

An important part of the design of any pressure-relief-valve design is its 
so-called {\em blowdown effect} \cite{Song2013,2010s:7,Safrany} in which the valve opens at a given {\em set pressure} $p_{\rm set}$ but only closes
when the pressure falls to a lower closing pressure, say
5 or 10\% below $p_{\rm set}$. Such a hysteresis loop between valve opening and closing avoids rapid cycling, in which a gradual build-up of pressure would cause the valve to open, the over-pressure to be rapidly relieved, and the valve to close almost instantaneously, only for pressure to build up again.  Indeed  even with blowdown, rapid cycling can occur if there is velocity-dependent pressure loss as fluid flows in the inlet piping to the valve. Such a possibility has led to the aforementioned 3\% rule, to ensure the line pressure loss is covered by the blowdown.

Another key part of the design of a DSOPRV is a limit to the maximum {\em lift} $x$ of the valve, usually implemented via an upper stop.  
This defines
{\em full lift} $x=x_{\rm max}$.
When choosing a valve for a particular application, it is important to understand its {\em capacity} $\dot{m}_{\rm in,cap}$ defined
as the mass flow rate through the valve 
at full lift and  10\% over the set pressure, that is  $1.1p_{\rm set}$. 

Typically, when specifying valve requirements for a particular plant, care is taken to choose the optimal set pressure and capacity for the process at hand. Most valves are designed to open by different amounts $x$, in proportion to the over-pressure $p-p_{\rm set}$, reaching full lift $x=x_{\rm full}$ at  $p-p_{\rm set}=0.1$. 

\begin{figure}[t]
	\centering
	\subfloat[]{
		\includegraphics[width = 0.25 \linewidth]{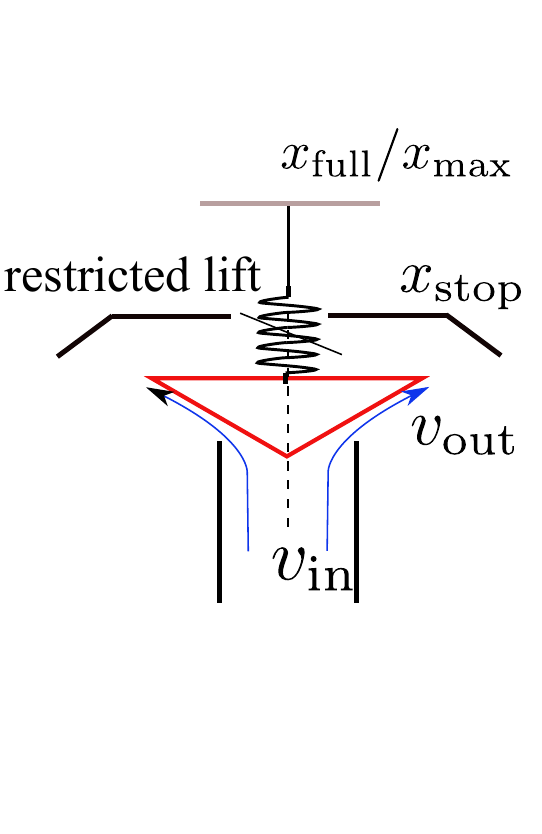}
	}
	\subfloat[]{
		\includegraphics[width = 0.6 \linewidth]{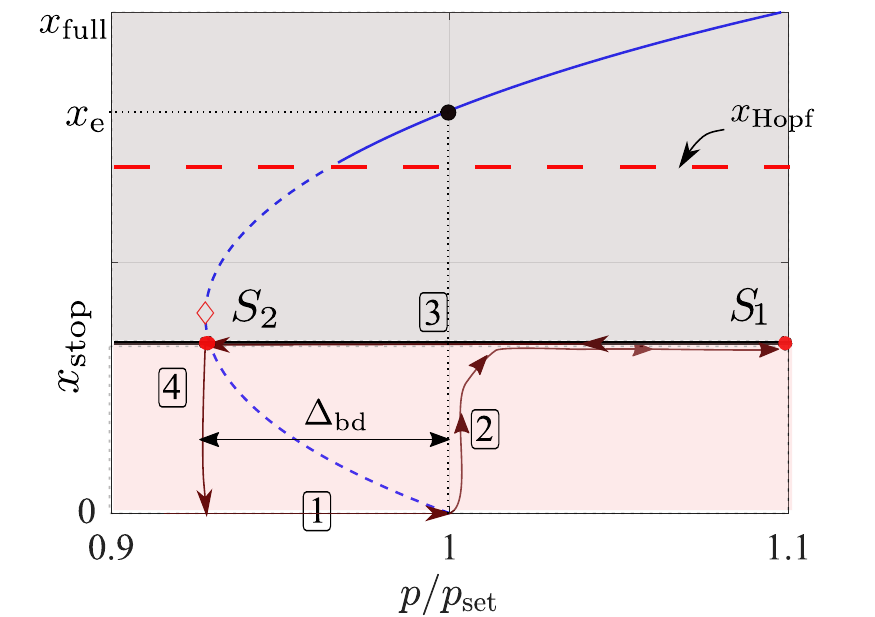}
	}
	\caption{
		The mechanistic working mode for an oversized valve with blowdown effect and upper stop:
		{\tikzsymbol[rectangle]{minimum width=6pt,minimum height=6pt,fill=gray}} virtual lift region;
		{\tikzsymbol[rectangle]{minimum width=6pt,minimum height=6pt,fill=pink}}  admissible region but dynamically unstable;
		{\protect\tikz \protect\draw [line width =
			1pt,color = blue] (0,0.3ex) -- (0.5cm,0.3ex);} stable equilibrium branch; {\protect\tikz \protect\draw [line
			width =
			1pt,dashed, color = blue] (0,0.3ex) -- (0.5cm,0.3ex);} unstable equilibrium branch;
		{\protect\tikz \protect\draw [line width =
			1pt,color = black] (0,0.3ex) -- (0.5cm,0.3ex);} location of the upper stop;
		{\protect\tikz \protect\draw [line width =
			1pt,color = red, dashed] (0,0.3ex) -- (0.5cm,0.3ex);} critical valve lift, such that without the upper stop, the valve would be unstable for all lift values below this line;
		{\color{red}{$\diamondsuit$}} the fold point;
		$\rightarrow$ ideal path segments of a pressure-relief event. More details including parameter values used to produce this simulation are given in \cref{sec:stability}.
	}
	\label{fig:stopper_below_fold}
\end{figure}

Nevertheless, as argued in \cite{HosCh2017-IV}, for most valves, the curve of equilibrium valve lift against pressure must be non-monotone. In particular, in order to allow blowdown to occur, there must be a static instability of the equilibrium lift as a function of pressure.  Specifically,  on quasi-static reduction of pressure through set pressure, there must be a stable branch of equilibrium valve lift which subsequently loses stability in a fold bifurcation (tipping point) leading to the rapid blowdown and full valve closure \cite{BazsoHos2015}.

The central idea developed in this paper is to make this blowdown effect more extreme, so that the valve would always rapidly pop open to full lift $x=x_{\rm max}$ as pressure is increased through p$_{\rm set}$. As we shall see, such behavior can be achieved by designing an appropriate
$A_{\rm eff}(x)$, which in turn is influenced through the effective angle of discharge from the valve. 
The idea illustrated in Fig.~\ref{fig:stopper_below_fold} is partially inspired by the results of
Smith and Desai in \cite{Smith2025}, who showed in a wide range of laboratory experiments that stability is enhanced by 
restricting the valve lift.  In Fig.~\ref{fig:stopper_below_fold} 
the maximum lift is restricted by lowering the upper stop so that it becomes a fraction of the full lift (about 40\% in the case illustrated, that is $x_{\rm stop} = 0.4x_{\rm full}$). 

Note from the figure that all stable equilibrium valve positions now occur for unreachable lift values, beyond the restriction. 
The requirement of most safety standards that the valve must open unimpeded means that the opening and closing of the valve is typically undamped. However, there is no requirement preventing energy absorption being added at the upper stop. As is common in impact mechanics \cite{Bro16}, if we model the contact with the upper stop via a simple restitution law, we are able to show in what follows that provided the restitution coefficient isn't too large, then the position in which the valve is pushed against the restriction via the fluid flow represents a stable equilibrium. Here, using the language of piecewise-smooth dynamical systems theory \cite{DiBernardo2007} we say that such a state is a stable {\em pseudo equilibrium}. This stability is achieved even if the equivalent unrestricted valve were to be in the region that is unstable to the quarter-wave instability (below the red dashed line in Fig.~\ref{fig:stopper_below_fold}).

To understand the diagram, note the two special pseudo-equilibrium points;
$S_1: (x,p)=(x_{\rm stop},1.1 p_{\rm set})$ which corresponds to the same valve capacity as in the standard design, and $S_2: (x,p) =(x_{\rm stop}, p_{\rm fold})$, which is the point where blowdown would occur. 
Thus a typical pressure-relief event would typically proceed via a {\em hysteresis loop} following the four stages illustrated by the corresponding numbers in Fig.~\ref{fig:stopper_below_fold}:
\begin{enumerate}
	\item The valve is closed but tank pressure accumulates until it reaches the set pressure;
	\item The valve then quickly pops open and settles on the upper
	stop via a brief low-energy chattering sequence;
	\item As the pressure continues to increase, the valve remains at the upper stop, along the pseudo equilibrium branch $S_1$--$S_2$. At some point, before reaching $S_2$, the pressure starts to decrease, while the valve remains in pseudo equilibrium;
	\item  The pressure decreases through $p_{\rm set}$ and continues towards the point $S_2$. At this point, the pseudo equilibrium becomes unstable and the valve now quickly closes. As the pressure is now below $p_{\rm set}$, we say that the over-pressure has been relieved. 
\end{enumerate}
We shall return to a full analysis of the dynamics in Fig.~\ref{fig:stopper_below_fold}
in \cref{sec:stability} below. 

\subsection{Outline} 

The rest of the paper is organized as follows. First, in \cref{sec:valve} we introduce the fundamental physics of DSOPRVs, including  derivation of a coupled mathematical model of a
DSOPRV connected to an inlet pipe. The resulting coupled valve and pipeline dynamics model (PDM) corrects several anomalies found in the equations presented in earlier studies \cite{HosCh2014-I,HosCh2015-II,HosCh2016-III,HosCh2017-IV}. We then show through CFD how complex discharge conditions can be captured via such a modeling approach through the introduction of a single discharge angle. This in turn leads to the effective-area-versus-lift $A_{\rm eff}(x)$ characteristics that underlie the hysteresis associated with blowdown. Section \ref{sec:eq-and-statbility} then recalls how to find the valve equilibrium lift for a given flow rate and analyze its stability. In particular, we reproduce the quarter-wave model (QWM) \cite{HosCh2015-II} that captures the dynamic instability mechanism that affects DSOPRVs at low flow rates.
We again correct anomalies in previous versions to more accurately model inlet pressure loss. This model is shown to accurately predict instability thresholds computed by the PDM model. 
\Cref{sec:stability} then analyzes the proposed stability enhancement in detail. Simulation results using the PDM show that for sufficiently small restitution coefficient of the upper stop, the stability threshold for a given flow rate is significantly enhanced. 
Finally,  \cref{sec:discussion} draws conclusions and suggests avenues for further research. 


\section{Mathematical modeling} 
\label{sec:valve}

The purpose of this section is to set out the 
basic physics of a typical DSOPRV,
and to explain from first principles how to derive its equations
of motion, correcting some errors or misconceptions in the literature, including in our earlier papers
\cite{HosCh2014-I,HosCh2015-II,HosCh2016-III,HosCh2017-IV}. 
We note that the valve should not be considered in isolation; by its nature, its dynamics is intimately coupled to 
the dynamics of the upstream fluid. 
For simplicity, we do not consider here any outlet piping, assuming that the fluid is vented into ambient conditions. We note however that in some situations, see e.g.~\cite{Chabane2009,Schmidt_outlet}, downstream effects can be important.

\begin{figure}
	\centering
	\subfloat[]{
		\includegraphics[width = 0.5 \linewidth]{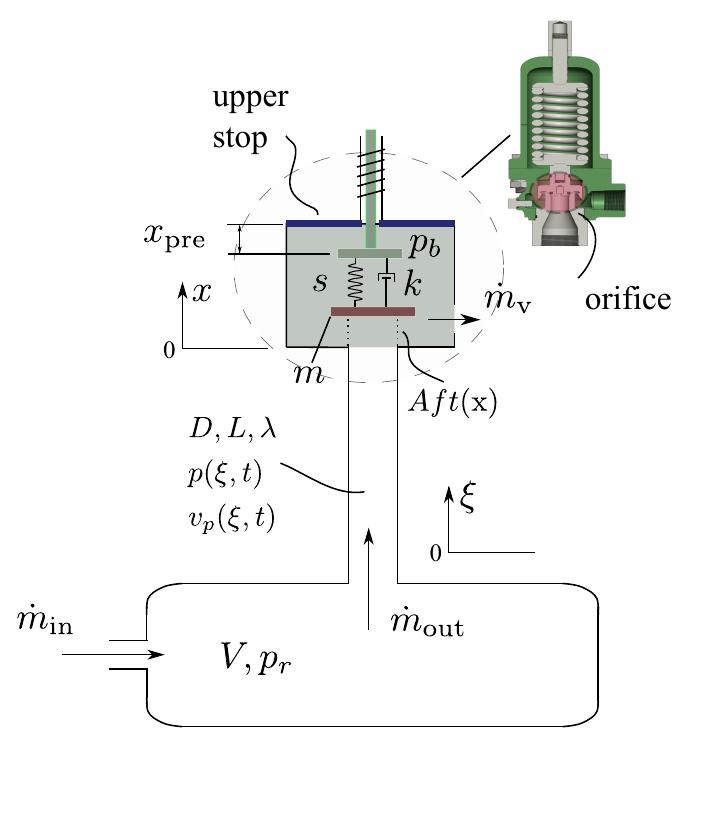}
	}
	\subfloat[]{
		\includegraphics[width = 0.4 \linewidth]{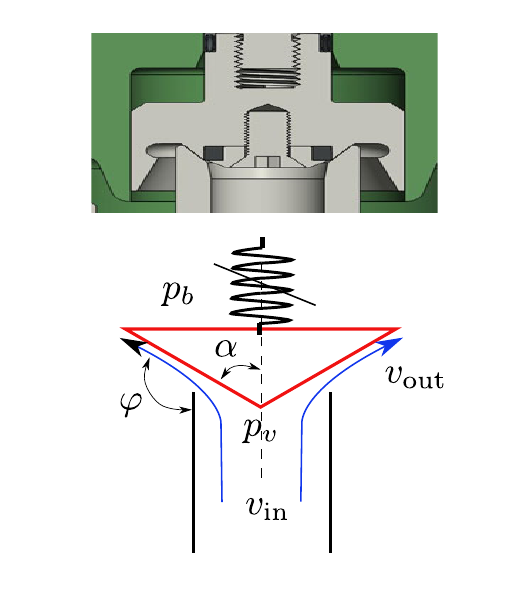}
	}
	\caption{(a) Definition sketch of pressure-relief system, with zoom at the valve in (b); see text for details. Insets show physical diagrams, taken from \cite{valve-inset}.}
	\label{fig:PSV-pipeline-sketch}
\end{figure}

We consider the generalized geometry of a pressure-relief system comprising a DSOPRV connected via a straight, rigid inlet pipe of length $L$ to a reservoir tank of volume $V$; see Fig.~\ref{fig:PSV-pipeline-sketch}. The meaning of the 
dynamic variables of the model are set out in 
Table~\ref{tab:variables} and the corresponding physical parameters are specified in Table~\ref{tab:parameters}. All data are consistent with a so-called 2J3 valve under either liquid or gas service.

The following subsections derive equations of motion for the three main components of the system --- the valve motion,
pipeline fluid dynamics, and tank pressure dynamics --- as well as provide appropriate boundary conditions to couple each component.

\begin{table}
	\begin{center}
		\begin{tabular}{|l|l|p{4in}|}
			\hline 
			Variable & Units & Meaning \\
			\hline 
			$t$ & sec & time \\
			$\xi$ & m & distance along pipe: $\xi=0$ at reservoir; $\xi=L$ at valve \\ 
			$x(t)$ & m & valve body displacement (lift): $0 \leq x < x_{\rm max}$ \\
			$y(t)$ & [-] & dimensionless valve lift $y=4x/D$ \\
			$\dot{m}_{\rm in}(t) $ & $\rm{kg}\:\rm{m}^3\rm{sec}^{-1}$ 
			& mass flow rate into reservoir \\ 
			$\dot{m}_{\rm out}(t)$ & $\rm{kg}\:\rm{m}^3\rm{sec}^{-1}$ 
			& output mass flow rate from reservoir into pipe \\
			$p_r(t)$ & Pa & reservoir pressure \\
			$\dot{m}_{v}(t)$ & $\rm{kg}\:\rm{m}^3\rm{sec}^{-1}$ & mass flow rate through the valve\\
			$p_v(t)$ & Pa & pressure at the valve inlet  \\
			$p(\xi,t)$ & Pa & inlet pipe fluid pressure \\
			$v(\xi,t)$ & $\rm{m}\:\rm{sec}^{-1}$ & inlet pipe fluid velocity
			\\ \hline 
		\end{tabular}
		\caption{Dynamic variables used in the mathematical model.  }
		\label{tab:variables}
	\end{center}
\end{table}

\begin{table}
	\begin{center}
		\begin{tabular}{|l|l|l|p{3.5in}|}
			\hline 
			Parameter & Units & Value & Meaning \\
			\hline 
			$\dot{m}_{in}$ & $\mathrm{kg/s}$ & 0-1.1847 & mass flow rate  \\
			$k$ & $\mathrm{N/m}$ & $5000$ & valve spring constant \\
			$c$ & $\mathrm{N s/m}$& $20$ & valve damping constant \\
			$m$ & $\mathrm{kg}$ & $0.45$ & total mass of the valve body (its moving part)  \\
			$x_{\rm pre}$ & ${\rm m}$ & $p_{\rm set} A/k$ & valve spring pre-compression \\
			$x_{\rm max}$ & ${\rm m}$ & ${D/4}$ &  valve full lift \\
			$x_{\rm stop}$ & ${\rm m}$ &  &  location of upper stop; $0<x_{\rm stop} < x_{\rm max}$ \\
			$e_{0,1}$ & [-] & $0.2-0.8$ & restitution coefficient at valve seat and upper stop \\
			$L$ & $\mathrm m$ & $ > 0$ & inlet pipe length \\
			$D$ & $\mathrm{mm}$ & $32.05$ & inlet pipe diameter \\
			$\lambda$ & [-]& $0$ & coefficient of wall friction  per unit pipe length \\
			$\alpha$ & $\rm{rad}$ & $(0,\pi)$& half-cone angle of valve body \\
			$C_d$ & [-]& 0.93 & discharge coefficient of fluid exiting valve \\
			$C_{\kappa}$ & [-] &  & choking factor, see Eq.~\eqref{eq:C_kappa} \\
			$\varphi$ & $\rm{rad}$ & $(0,\pi)$ & discharge angle of fluid from the valve\\
			$p_b$ & $\rm bar$ & 1 & ambient pressure  \\
			$p_{\rm set}$ & $\rm bar$ & 5 & ambient pressure  \\
			$T_{0}$ & $\rm K$ & 293 & ambient temperature  \\
			$R$ & $\rm{J/(kg K)}$ & 288 & gas constant \\
			$\kappa$ & [-] & 1.4 & special heat ratio \\
			$\rho$ & $\rm{ kg/m^3}$&  & (near) constant density of fluid \\
			$a$ & $\rm m/s$ & & sonic velocity of fluid, see Eq.~\eqref{eq:sonic}\\
			$V$ & $\rm m^3$ & 1 & reservoir volume  \\
			\hline 
		\end{tabular}
		\caption{Meaning of the physical parameters in the mathematical model, together with their default values used in this paper, unless otherwise stated
		}
		\label{tab:parameters}
	\end{center}
\end{table}

\subsection{Valve dynamics}

The main moving part of a DSOPRV can be modeled as a single degree-of-freedom mechanical 
system 
\begin{equation} \label{eq:valve-dynamics-equations} 
	m \ddot x + c \dot x + k  ( x_{\rm pre}  + x)  = F_{\rm fluid}    
	, \qquad \text { for } 0 < x < x_{\rm stop}, 
\end{equation}
where a dot means differentiation with respect to time and $F_{\rm fluid}$ represents all fluid forces acting on the valve body. We assume that the valve seat and the upper stop are rigid and the dynamics there can be described by an 
instantaneous Newtonian restitution law. Hence we assume that the valve body velocity $\dot{x}^+$ immediately following an upper or lower impact can be written in terms of its pre-impact velocity $\dot{x}^-$ via 
\begin{align}
	\dot{x}^+ &= - e_0 \dot{x}^-, \quad \mbox{for } \dot{x}^- <0 \mbox{ and } x=0, \label{eq:lowerImpact} \\
	\dot{x}^+ &= - e_1 \dot{x}^-, \quad \mbox{for } \dot{x}^- >0 \mbox{ and } x=x_{\rm stop}. \label{eq:upperImpact}
\end{align}

The fluid forces can be decomposed into a static part (which applies equally if the valve is open or closed) and a dynamic part arising from the change in momentum of the fluid velocity passing through the valve. Hence we can write
\begin{align}
	F_{\rm fluid} &\coloneqq  (p_v - p_b) A_0 + \dot{m}_v(v_{\rm in} + v_{\rm out} \cos\varphi).
	\label{eq:Ffluid}
\end{align}
where 
$$
p_v \coloneqq p(L,t), \qquad A_0 = \pi D^2/4
$$
are the static pressure of the fluid on the valve and the valve seat area, respectively $\dot{m}_v(t)$  is the mass 
flow rate through the valve, $v_{\rm in} = v(L,t)$ and $v_{\rm out}(t)$ is the magnitude of outflow velocity from the valve. (Note that the discharge angle is defined with respect to the downwards positive direction, so that if the valve were absent $\cos(\varphi) =-1$ and there would be no dynamic component to the fluid force). 

The measurement of $v_{\rm out}$ is unfortunately difficult in practice, and the prediction of $v_{\rm out} (t)$ from the other dynamical variables typically depends on many features of the valve geometry; see~\ref{sec:Appendix_I}.
An alternative approach, used in our earlier work 
\cite{HosCh2017-IV,Hos2014}, was inspired by the industrial practice that 
force on the valve and the pressure at the valve inlet are amenable to experimental measurement. Thus a convenient way to characterize the fluid force is via the definition of an {\em effective area} 
\begin{equation} \label{eq:Aeff-estimation}
	F_{\rm fluid} = (p_v-p_b) A_{\rm eff}(x).
\end{equation}
Here $A_{\rm eff}(x)$, although having dimensions of area, is not a physical area, but rather can be thought of as a curve that defines the valve's operational characteristics. In truth, $A_{\rm eff}$
is a function of other variables such as inlet velocity, but we make the assumption that for a given fluid, for each valve lift $x$ there is a unique equilibrium of the flow through the valve (be it stable or otherwise) such that we can write $A_{\rm eff}$(x). We shall explore this approach further in \cref{sec:effective_area}, relating $A_{\rm eff}$ to more familiar concepts such as discharge coefficient $C_d$ and angle $\varphi$.  We shall also see in  \cref{sec:eq-and-statbility} that the effective area characteristics provide a useful way of explaining the complex flow physics that underlie the blowdown effect. 

Note from the form of Eq.~\eqref{eq:valve-dynamics-equations} that
the set pressure is easily calculated to be
\begin{equation}
	p_{\rm set} = p_b + \frac{ kx_{\rm pre}}{A_{0}} \label{eq:pset_def}
\end{equation}
%

\subsection{Inlet pipe dynamics} 

For simplicity, we shall ignore effects of gravity and assume
that the inlet pipe lies along the $\mathbf{j}$ axis. Let $\xi\in [0,L]$ represent a co-ordinate along the pipe axis, with $\xi=0$  at the tank end and $\xi=L$ at the valve. 
Following \cite{HosCh2015-II}, assuming that the pipe is sufficiently slender, we use a 1D fluid dynamics approximation and denote  by $p(\xi,t)$ and $v(\xi,t)$ the pressure and velocity distribution respectively along the pipe axis. 

If we restrict attention to sufficiently small variations to steady flows within the pipe, so that we can assume near constant density and temperature, then the equations for conservation of mass, energy and momentum can be reduced to the standard equations of 1D pipe-flow
\cite{HosCh2015-II}:
\begin{subequations}
	\label{eq:pipeline-dynamics-equation}
	\begin{align}
		\frac { \partial p } { \partial t } + v \frac { \partial p } { \partial \xi } + a^{ 2 } \rho \frac { \partial
			v} { \partial \xi } & = 0,
		\\
		\frac { \partial v } { \partial t } + v \frac { \partial v } { \partial \xi }  +  \frac { 1 } { \rho } \frac {
			\partial p } { \partial \xi } +  \frac {\lambda } { 2D } v |v| & = 0. 
	\end{align}
\end{subequations}
Here $D$ is the inlet pipe diameter, $\rho$ is the assumed (near) constant density of the fluid, $a$ its sonic velocity which satisfies
\begin{equation}
	a^2 =
	\begin{cases}
		\kappa R T \quad \mbox{for gas,}\\
		E/\rho \quad \mbox{for liquid,}
	\end{cases}
	\label{eq:sonic}
\end{equation}
where $\kappa$ and  $R$ are the usual Boltzmann and real gas law constants, $T$ is
temperature (which is assumed to be constant) and $E$ is the bulk modulus. Also, 
$\lambda$ is a coefficient of wall friction, 
(see \cite[eq.~(3.63)]{Zucker2002}).
Note that here we have corrected a misprint in \cite{Hos2014,HosCh2015-II} which had the incorrect sign of the wall friction term. 

\subsection{Tank pressure dynamics and boundary conditions}

The tank pressure dynamics can be assumed to 
satisfy a simple mass balance equation
\begin{equation} \label{eq:tank-pressure-dynamics}
	V \dot{p}_r =  a^2 (\dot{m}_{\rm in} - \dot{m}_{\rm out}),
\end{equation}
where we think of $\dot{m}_{\rm in}(t)$ as the net mass inflow that causes the overpressure that leads to the valve opening. 
Conservation of mass shows that the outflow must be equal to the mass flow into the pipe:
\begin{equation} \label{eq:mass-flow-rate-of-tank}
	\dot{m}_{\rm out} = \frac{\pi D^2}{4} \rho v(0,t), 
\end{equation}
leading to a simple velocity boundary condition at the $\xi=0$
end of the pipe. 

Determination of pressure boundary conditions at $\xi=0$
needs more care, as it depends on the kind of fluid being serviced. 
For the case of ideal gas flow, which is compressible, we can use the principle of conservation of energy, assuming tank pressure to be approximately constant and the pipe flow to be isentropic. The  details are given in \cite[eq.(13)]{Hos2014}, and lead to an equation of the form 
\begin{equation}\label{eq:BC_1_at_0}
	p(0,t) = p_r(t)\left[ 1  - \frac{v(0,t)^2}{2 c_p T} \right]^{\frac{\kappa}{\kappa -1}} := p^{\rm (gas)}_0(t)
\end{equation}
where $c_p$ is the specific heat capacity of the fluid at fixed pressure (assumed to be $p=p_{\rm set}$). 
Alternatively, for liquid service valves we can assume incompressibility and apply Bernoulli's equation to give 
\begin{equation}\label{eq:BC_2_at_0}
	p(0,t) = p_r(t) - \frac{\rho}{2} (v(0,t))^2 \coloneqq  p^{\rm (liquid)}_0(t).
\end{equation}
In what follows we shall define the unified inlet boundary condition 
\begin{equation} \label{eq:BC_at_0}
	p(0,t) = p_r(t) \left [1- \chi(p_r,v(0,t)) \right ],
\end{equation}
where the {\em pipe inlet pressure loss} $\chi$ is defined by
Eq.~\eqref{eq:BC_1_at_0} or \eqref{eq:BC_2_at_0} depending on whether liquid or gas is being serviced. Note that typically 
$\chi \ll 1$.

At the valve end of the pipe, it is the pressure boundary condition that is trivial:
\begin{equation}
	p(L,t) = p_v(t) \label{eq:p_v}
\end{equation}
For the velocity boundary condition, however, we need to use conservation of mass, by
defining the {\em flow through area} $A_{\rm ft}(x)$ of the valve, which is the cross-sectional area of the  narrowest part of the valve through which fluid flows when the valve is at lift $x$.  

In general, an expression for $A_{\rm ft}$ needs to be determined from
the detailed geometry of the valve. However for the simple case of a conical valve, then using the half-cone angle $\alpha=\pi-\varphi$ as in Fig.~\ref{fig:PSV-pipeline-sketch}(b) we find from simple geometry that 
\begin{equation}
	\label{eq:cone-shape-Aft}
	A_{\rm ft}(x) \coloneqq \pi x \sin \varphi (D - x \cos\varphi \sin\varphi). 
\end{equation}
Note that for the typical case  $x \ll D$ we see that $A_{\rm ft}(x) \propto x$ for small valve lift $x$.

Then conservation of mass flow at the valve end of the pipe gives 
$$
A_0 v(L,t) = A_{\rm ft}(x) v_{\rm out},
$$
where $v_{\rm out}$ is the outflow from the valve which can be calculated from first principles using Bernoulli's principle (see \ref{sec:Appendix_I} for details). 
We obtain
\begin{equation}\label{eq:BC-V-at-1}
	v(L,t) \coloneqq v_L(t) =
	\frac{A_{\rm ft}(x)}{A_0} v_{\rm out} = 
	C_d C_\kappa \frac{A_{\rm ft}(x)}{A_0} \sqrt{ \frac{p^{*} }{\rho}},
\end{equation}
where 
$$
p^*: =\begin{cases} 
	p_v & \mbox{for gas}, \\
	p_v-p_b & \mbox{for liquid.}
\end{cases}
$$
Also, $C_d$ is the commonly used discharge coefficient that measures orifice losses through the valve exit 
and $C_\kappa$ is a geometric term that depends on whether the flow through the valve is choked (for the case of gas) or pressure-driven (for liquids).  
\begin{equation} 
	C_{\kappa} = \begin{cases} 
		\: 	\displaystyle\sqrt{ \kappa \left (\frac{2}{1+\kappa}\right ) ^{\frac{\kappa+1}{\kappa-1}}}  & 
		\mbox{for gas,}  \\
		\: \sqrt{2}, &  
		\mbox{for liquid.}
	\end{cases}
	\label{eq:C_kappa}
\end{equation}
Here $\kappa = c_p/c_v$ is the ratio of specific heat capacities at constant pressure and volume. The corresponding mass flow through the valve is
\begin{equation}
	\dot{m}_v = \rho A_0 v_L.
	\label{eq:mv}
\end{equation}

\subsection{Simulation method}

Given expressions for $A_{\rm eff}(x)$ and $A_{\rm ft}(x)$ that
represent a specific valve geometry, and a prescribed 
$\dot{m}_{in}(t)$, then the Eqs.~\eqref{eq:valve-dynamics-equations},
\eqref{eq:pipeline-dynamics-equation}
and \eqref{eq:tank-pressure-dynamics}
represent a
formally well-posed system of ordinary and partial differential
equations (ODEs/PDEs) for the dynamical variables $x(t)$, $p_r(t)$,
$v(\xi,t)$ and $p(\xi,t)$ coupled via the boundary conditions Eqs.~\eqref{eq:BC_at_0} and \eqref{eq:BC-V-at-1}, with additional impact conditions Eqs.~\eqref{eq:lowerImpact} and \eqref{eq:upperImpact}.
Having solved these equations, the mass flow variables   $\dot{m}_{\rm out}$, $\dot{m}_v$ and $p_v$ are defined in terms of these unknowns via the expressions Eqs.~\eqref{eq:mass-flow-rate-of-tank}, 
\eqref{eq:mv} and \eqref{eq:p_v}. 
We shall refer to this model in what follows as the full
{\em valve-and-pipeline dynamics model} (PDM).  

To simulate the system, the ODEs describing the reservoir and valve dynamics are solved using a fourth-order Runge-Kutta scheme with error-controlled adaptive stepsize. The PDE part is solved using different approaches for the gas and liquid cases.  For liquids, we use the method of characteristics \cite{wiley1982fluid} to track the pressure and velocity changes inside the pipe. For gases, a standard two-step Lax–Wendroff method 
\cite{Zucker2002} is used instead.
Further details of the implementations are given in \cite{HosCh2014-I} 
and \cite{HosCh2016-III} respectively, where experimental validation is provided  for both liquid and gas service valves.

A numerical code has been constructed in \verb*|C++| to solve the PDM using these methods. 
It is available at 
\hyperref{https://github.com/HongTang973/PS-simulator}{}{}{https://github.com/HongTang973/PS-simulator}.
In particular, we have updated the numerical toolbox developed in
\cite{Hos2014} which was already validated by comparison with experimental data in \cite{HosCh2014-I,HosCh2016-III}.

In the new implementation, it was found necessary to 
resolve impact events Eqs.~\eqref{eq:lowerImpact},\eqref{eq:upperImpact} more accurately to avoid spurious oscillatory solutions when the valve is close to its upper stop.
To do this, we used 
event detection using nonlinear interpolation to accurately locate the
events, and the time stepping is rolled back and restarted to carry on integration.
Also, contact with either the upper or lower stop 
is not achieved instantaneously, but via a decaying  
chattering sequence of impacts, with 
small impact velocities $|\dot{x}(t^-)|$ at $x=x_{\rm stop}$ or $x=0$. To avoid problems with trying to resolve infinitely many impacts in finite time, the algorithm introduced in \cite{Nordmark_geometric}
was used to jump straight to contact at the convergence point of
the geometric sequence once the impact velocity becomes too small. 
This contact is then retained for future timesteps until $\ddot{x}$ changes sign. These extra ingredients are necessary in order to accurately resolve pseudo equilibrium solutions where the valve body is permanently in contact
with the upper stop.

\section{Valve geometry and effective area curve}
\label{sec:effective_area}

The purpose of this section is to motivate how the concept of the effective-area-versus-lift curve can capture the complex flow physics within different valve geometries. 

One way of representing the effect of an open valve is via
its discharge coefficient
$C_d$ and discharge angle $\varphi$. From these two quantities, it is possible to define an analytic
expression for $A_{\rm eff}(x)$. The details of the derivation
are explained in 
\ref{sec:Appendix_I}, where we obtain 
\begin{equation}
	A_{\rm eff}(x) =  \begin{cases} 
		A_0+ C_d^2 C_{\kappa}^2 \displaystyle\frac{ A_{\rm ft}(x)}{A_{\rm 0}} (A_{\rm ft}(x) +
		A_{0} \cos\varphi),
		& x>0 \\
		A_{0}, & x=0 \end{cases}.
	\label{eq:Aeff_from_Aft}
\end{equation}
This expression depends on $x$ through the expression for $A_{\rm ft}(x)$ which
we recall from Eq.~\eqref{eq:cone-shape-Aft} can be written as a linear function of
$x$ plus a quadratic correction for such simple geometries.
Hence $A_{\rm eff}(x)$ can be written to leading order as a quartic polynomial in dimensionless lift $y = 4x/D$:
\begin{equation}
	A_{\rm eff}(y) = A_0 + \tilde a_1 y + \tilde a_2 y^2 + \tilde  a_3 y^3 + \tilde  a_4 y^4,
	\label{eq:Aeff_Taylor}
\end{equation}
where expressions for the coefficients $\tilde a_1,\ldots, \tilde a_4$ in the case of the discharge angle $\varphi$ are given in the \ref{sec:Appendix_I}.

It should be emphasized, however, that estimating the
effective-area-versus-lift characteristics (and discharge
coefficients) is the cornerstone of valve modeling, and such
analytical approximations rarely suffice for accurate valve
calibration in industrial settings, bearing in mind the detailed
turbulent flow physics.  Instead, practitioners typically resort to
empirical measurements, backed up by CFD.  For the computational part,
Reynolds-averaged Navier-Stokes (RANS) equations are commonly used,
because they provide a sufficient balance between accuracy and
computational times; see for example
\cite{Chabane2009,Song2013,Song2014,Wu2015}. We have
followed that approach here in order to test the effectiveness of our
analytical approximation of $A_{\rm eff}(x)$.

Specifically, we have used commercial CFD software to 
compute $A_\mathrm{eff}(x)$ and $C_d$ for some prototypical valve shapes,
using a quasi-2D approach that exploits the axial symmetry of the geometries.
The details of the method are given in \ref{sec:Appendix_II}.
We have taken three disc-shaped valves with different jet angles, and one conical valve.
For reference, the inner pipe diameter was taken to be $D=40.2$ \, mm in all cases.

\begin{figure}
	\vspace{-4em}
	\centering 
	\includegraphics[width=1.0\textwidth]{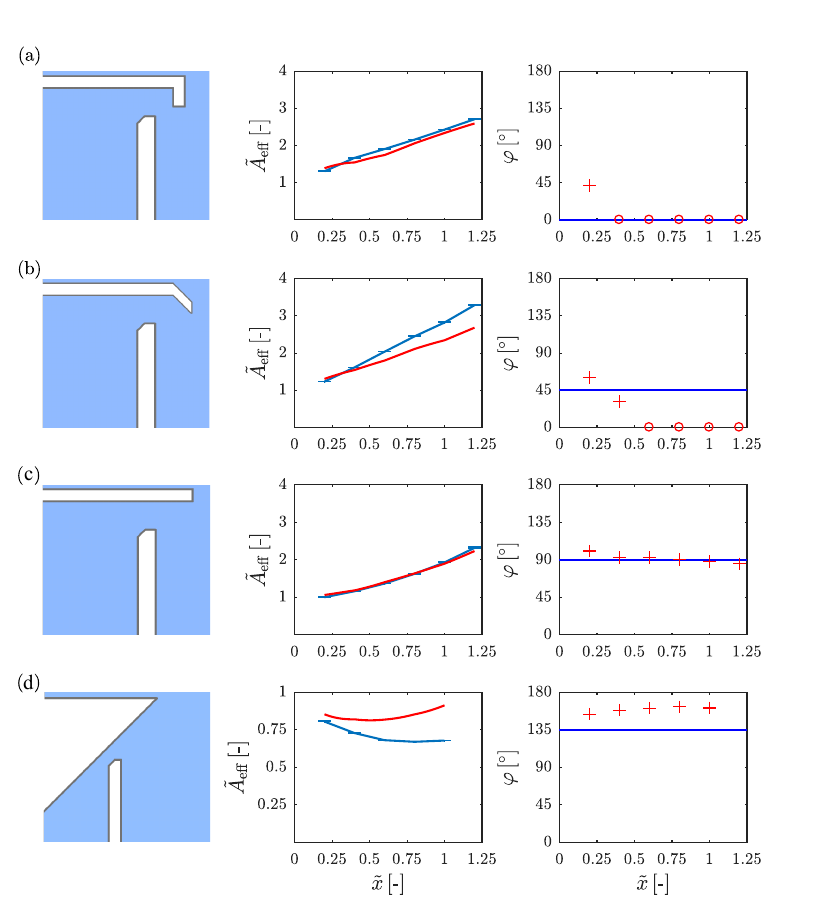}
	\caption{Effective-area characteristics and fitted discharge angles for different valve geometries. In each row, the left panel shows the valve geometry, the middle panel shows the normalized effective area $\hat{A}_{\rm eff}=A_{\rm eff}/A_0$ as a function of normalized lift $\tilde{x}=x/x_{\max}$, and the right panel shows the corresponding fitted discharge angle $\varphi$. In the middle panels, CFD results are compared with the analytical estimate. In the right panels, markers (+) denote the fitted values of $\varphi$, while horizontal lines indicate the geometrical deflection angles. Where fitting was not possible, the closest approximation (o) is shown.}
	\label{figAeffJetAngleGrid}
\end{figure}

The second column of Fig.~\ref{figAeffJetAngleGrid} depicts the computed effective area
curves with their analytical estimation based on
Eq.~\eqref{eq:Aeff_Taylor} and data in 
Table \ref{tab:parameters}.  We found that the effective area does not
depend significantly on the pressure difference; the error bars in the plots represent the range in which the effective areas varied for a
fixed lift value but for different inlet pressures. For the $0^\circ$
and the $90^\circ$ cases, the match between the CFD and analytic
expression is exceptionally good. However, in the $45^\circ$ and the
$135^\circ$ cases there are significant differences. 

We have also
computed the discharge coefficient $C_d$ for these geometries. We find,  
as is well known in practice, that
a constant value of $C_d$ is not strictly valid and  typically is a decreasing function of discharge velocity. This result is not particularly important for the arguments presented here, so we relegate the details to
\ref{sec:Appendix_II}. 

The equivalent discharge angles $\varphi$ were calculated from the CFD
by solving Eq.~\eqref{eq:Aeff_from_Aft} for $\cos(\varphi)$, given the computed values of $A_{\rm eff}$ and $C_d$.
The third column of Fig.~\ref{figAeffJetAngleGrid} shows these computed angles  
(as red `+' markers)
compared to the expected angles from the geometry
(continuous blue lines). In some cases fitting was not possible.
For example, in the topmost case the geometrical deflection angle is
$0^\circ$, which gives the highest $\hat{A}_{\rm eff}$ for a given
$C_d$ according to Eq.~\eqref{eq:Aeff-estimation}, however, the
simulated $\hat{A}_{\rm eff}$ values are (marginally) higher than
their analytically estimated counterparts in all but one point,
therefore no valid $\varphi$ exists in these cases. For these points,
the closest approximations of $\varphi=0$ are shown in
Fig.~\ref{figAeffJetAngleGrid} (marked by red circles). Interestingly,
these differences are quite significant in the second ($\varphi =
45^\circ$) case --- the simulated $\hat{A}_{\rm eff}$ values are even
higher than in the $\varphi = 0^\circ$ case. This observation
suggests, given that the discharge coefficient is well-defined, there appears to exist an additional fluid-mechanical effect related to the
precise valve geometry. Nevertheless, the results confirm that
Eqs.~\eqref{eq:Aeff-estimation} and \eqref{eq:hat_Aeff} sufficiently
capture the effect of momentum forces and give a simple, yet
relatively accurate estimation for the effective area.

For valve geometries that are more complex than these simple disk or cone valves (e.g.~if a housing is present) more irregular effective area-versus-lift functions can be found, see e.g.~\cite{Song2013} and \cite{BazsoHos2015}. 
In these cases, we advise that either specific CFD simulations or
experimental measurements would still need to be conducted in order to find accurate approximations to 
$A_{\rm eff}(x)$.


\section{Valve equilibrium and instability}
\label{sec:eq-and-statbility}

This section investigates existence and stability of different equilibrium states of the valve model developed in Sec.~\ref{sec:valve}, in
the case of steady inflow $\dot{m}_{\rm in} ={\rm const}$.
There are three different kinds of possible equilibrium states :
\begin{enumerate}
	\item {\em Closed} $x \equiv 0$. In this state there is no flow;
	$v(\xi,t)=0$, $p(x,t)=p_r(t)$, where the tank pressure increases at the fixed rate $\dot{m}_{\rm in} /V$. Provided $p_r(t)< p_{\rm set}$ then this state is stable to small perturbations. 
	\item {\em Dynamic equilibrium} $x \equiv x_e$ where $0<x_{e}<x_{\rm stop}$. Here the flow in the pipe must be steady so that $p_r = \mbox{const.}$ and 
	$v(\xi,t) = v_s(\xi)$, $p(\xi,t)=p_s(\xi)$, where 
	$v_s(\xi)$ and $p_s(\xi)$ 
	are the steady time-independent solutions of Eq.~\eqref{eq:pipeline-dynamics-equation}.
	This equilibrium state need not necessarily be stable. 
	\item {\em Pseudo equilibrium}. Here $x \equiv x_{\rm stop}$ and the flow is
	again steady, as described in the previous case. In keeping with the definitions for piecewise smooth dynamical systems \cite{DiBernardo2007} we call this a pseudo equilibrium because it does not satisfy the
	Eq.~\eqref{eq:valve-dynamics-equations} with $\dot{x}=\ddot{x}=0$.
	Instead, there is a constraint force that keeps the valve in contact with its upper stop provided the fluid force $F_{\rm fluid} >  k  ( x_{\rm pre}  + x)$ so
	that $\ddot{x}$ described by Eq.~\eqref{eq:valve-dynamics-equations} remains
	positive. Again, a pseudo equilibrium need not necessarily be stable. 
	Note that a transition between an equilibrium and a pseudo equilibrium would occur if $\ddot{x}=0$, such a situation being known as a {\em boundary equilibrium bifurcation} \cite{DiBernardo2008,Hong2023,TangSimpson}.
\end{enumerate}

\subsection{Dynamic equilibrium}

We now investigate the properties of dynamic equilibria, and how this leads to an effective valve characteristic curve of lift versus tank pressure.
Static instability, or folds, of this curve are what underlies the blowdown effect. We also recall the fundamental quarter-wave dynamic instability found in previous works, 
and how the onset of this instability can occur either before or after the fold
with increasing valve lift. 

At dynamic equilibrium, the flow in the pipe is steady, and for simplicity let us suppose that pipe friction is negligible. (Inclusion of pipe friction
$\lambda \neq 0$ does not change the qualitative argument, it just adds an additional factor between the valve and reservoir pressures $p_v$ and $p_r$.)
For simplicity, we also assume that $p^*=p_v-p_b$ in Eq.~\eqref{eq:BC-V-at-1}, which would be true in the liquid 
case and is an approximation for gas service valves
(The true expressions in the gas case are qualitatively similar, but algebraically more cumbersome).
Then, from the boundary conditions, we have 
$$
p(L,t) = p(0,t):=p_v(t),
$$ 
$$
v(0,t)=v(L,t)= v_L =  C_d C_k\frac{A_{\rm ft}}{A_0} \sqrt{\frac{p_v-p_b}{\rho}}.
$$
The reservoir dynamics are therefore given by
\begin{align}
	\label{eq:3rd-system-ode}
	V \dot{p}_r & =   a^2 \left[\dot{m}_{\rm in} - C_d C_\kappa A_{\rm ft}(x) \sqrt{\rho (p_v-p_b)} \right].
\end{align}

A dynamic equilibrium can then be found by setting
$\dot{p}_r$, $\dot{x}$ and $\ddot{x}$ to zero in Eqs.~\eqref{eq:valve-dynamics-equations} and  \eqref{eq:3rd-system-ode}.
From the first of these equations we obtain
\begin{equation}
	p_v =  p_b+ \frac{ \dot m_{\rm in}^2}{\rho C_{\kappa}^2  C_d^2 A_{\rm ft}(x)^2 },
	\label{eq:equil1}
\end{equation}
and from the second 
\begin{equation}\label{eq:equil2}
	p_v  = p_b + \frac{1}{A_{\rm eff }(x)}[k ( x_{\rm pre}  + x)].
\end{equation}
Setting the right-hand sides of these two equations to be equal, we get a
nonlinear equation for $x$,  which can be solved to find the dynamic equilibrium $x=x_e$ for a given inflow parameter $\dot{m}_{in}$. Depending on the value of
$\dot{m}_{in}$, and the form of $A_{\rm eff}(x)$,
this equation may have a unique solution, multiple solutions, or no solution. 

To find the corresponding equilibrium value for $p_r$ we can use the relation
Eq.~\eqref{eq:BC_at_0} to write 
\begin{equation}
	p_r=\frac{p_v}{1 - \chi(p_r,v_L)}
	\label{eq:equil3}
\end{equation}
where $p_v$ is given by Eq.~\eqref{eq:equil1}. 
We note that Eq.~\eqref{eq:equil3}  leads to a well-defined solution for $p_r$, since  $\chi \ll 1$. 

Thus Eqs.~\eqref{eq:equil1}--\eqref{eq:equil3} can
effectively be used to define the equilibrium characteristic of the valve
and pipe system which can be expressed as the valve lift
$x_e$, as a function of $\dot{m}_{in}$,  $p_r$, or $p_v$. In 
simulations of a pressure-relief event $\dot{m}_{in}$ would be an input parameter,
and $p_r$, $p_v$ as dependent variables that have to be solved for. However, a valve manufacturer would more likely think of the valve characteristic as a curve of valve lift $x_e$ versus valve inlet pressure $p_v$. For the rest of this section though it is helpful to think of the valve characteristic as $x_e$ as a function of tank pressure $p_r$.

Note that by setting $p_r=p_{\rm set}$ and $v(\xi,t)=0$ we can find a solution
$x=0$. Then, using the implicit function theorem, provided that $A_{\rm eff}$ is a smooth function of $x$, it is straightforward to see that there
must be a single branch of dynamic equilibria $x_e$ that emerge from
$(p_r,x_e)$ into $x_e>0$. Also, because of the (approximately in the case of gas flow)  linear dependence of Eqs.~\eqref{eq:equil1}--\eqref{eq:equil3} on $p_r$
for each value of $x_e$ between $0$ and $x_{\rm stop}$, there will be a unique value of $p_r$, so we can write
\begin{equation}
	\label{eq:static_curve_pr_x_liquid}
	\mbox{valve equilibrium characteristic:} \quad  p_r = P(x_e)
\end{equation}
where the function $P$ is obtained from Eqs.~\eqref{eq:equil1}--\eqref{eq:equil3} upon
elimination of $\dot{m}_{\rm in}$ and $p_v$.  
Equation~\eqref{eq:static_curve_pr_x_liquid} defines the characteristic curve of \emph{reservoir pressure versus valve lift}.

\subsection{Static instability and blowdown effect} 
\label{sec:blowdown}
It is helpful to consider the stability of the equilibrium Eqs.~\eqref{eq:equil1} and \eqref{eq:equil2} under the assumption that the flow in the inlet pipe remains steady.
Such a {\em static instability} (sometimes also called a tipping point)
can be analyzed in the usual way through 
linearization of   Eqs.~\eqref{eq:valve-dynamics-equations} and \eqref{eq:3rd-system-ode} about 
$x=x_e$ and
$p_r=P(x_e)$, 
and demanding that eigenvalues of the associated Jacobian matrix should be in the left-half plane. 

However, there is a simplified approach to understand such 
\emph{static} instabilities just by looking at the form of the 
Eq.~\eqref{eq:valve-dynamics-equations} with 
$F_{\rm fluid} = (p_v - p_b) A_{\rm eff}$. Then 
$p_v$ can be written in terms of $p_r = P(x_e)$ using Eq.~\eqref{eq:equil3}, where we note that $p_v= p_r(1-\chi(p_r,v_L))$ where $\chi$ implicitly depends on $x$ through $P_r(x)$ and $v_L$, so with a slight abuse of notation we write $\chi=\chi(x)$. 

Next, linearizing 
Eq.~\eqref{eq:valve-dynamics-equations} 
about $x=x_e$ we get 
$$
m\ddot{x} + c \dot{x} = [P(x_e)(1-\chi) -p_b)A_{\rm eff} ^\prime]x
+ [P^\prime(1-\chi(x_e))-P(x_e)\chi^\prime]A_{\rm eff} x
- k x 
$$
where $\prime$ means differentiation with respect to $x$ at
$x=x_e$.  Thus we can write the linearized dynamics as
$$
m\ddot{x} + c\dot{x} + k_{\rm eff} =0,
$$
where the $k_{\rm eff}$ is the 
{\em effective stiffness} of the 
valve 
\begin{equation} \label{eq:effective-stiffness}
	k_{\rm eff} \coloneq
	k- P(x_e)(1-\chi) -p_b)A_{\rm eff}^\prime 
	- [P^\prime(1-\chi(x_e))-P(x_e)\chi^\prime]A_{\rm eff}
\end{equation}

In the industrially relevant case that the mechanical damping $c$ is negligible then a necessary condition for static stability is
$k_{\rm eff} >0$. In contrast if the effective stiffness is negative, then the dynamic equilibrium will be unstable. Thus, the condition for the onset of such a static instability will be given by 
$k_{\rm eff} =0$.

From a dynamical systems perspective such a condition defines
a {\em saddle-node bifurcation}, which can be identified by a fold point in the characteristic curve of $x_{e}$ as a function of $p_r$. Such an instability (labeled
$S_2$ in Fig.\ref{fig:stopper_below_fold}) when reached upon slowly reducing the pressure from a statically stable equilibrium, would lead to a sudden jump in the valve lift.
If this occurs at a pressure below set pressure, this would cause the valve to close.

\begin{table}
	\centering
	\caption{Effective-area coefficients in Eq.~\eqref{eq:Aeff_Taylor} for different half-open angles $\alpha$. 
		The first row shows the corresponding valve geometries.}
	\label{tab:-aeffs-from-different-open-angles}
	\begin{tabular}{l c c c c}
		\toprule
		& \multicolumn{3}{c}{Half-open angle $\alpha$} \\
		\cmidrule(lr){2-5}
		& $\pi/2$ & $\pi/3$ & $3\pi/5$ & $2\pi/3$ \\
		\midrule
		
		Geometry
		& \raisebox{-0.5\height}{\includegraphics[width=0.14\linewidth]{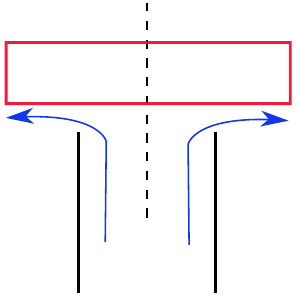}}
		& \raisebox{-0.5\height}{\includegraphics[width=0.14\linewidth]{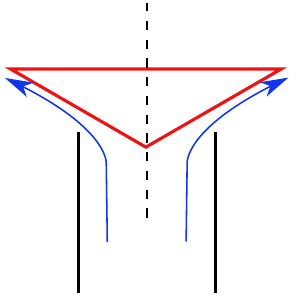}}
		& \raisebox{-0.5\height}{\includegraphics[width=0.14\linewidth]{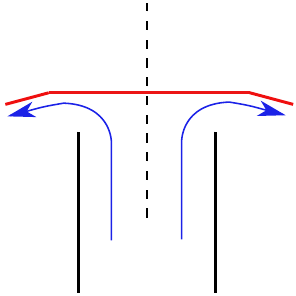}}
		& \raisebox{-0.5\height}{\includegraphics[width=0.14\linewidth]{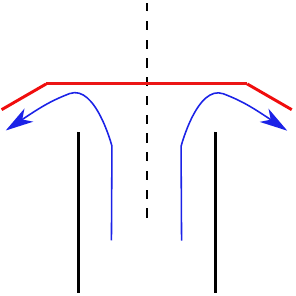}}
		\vspace{1em} \\
		
		$\tilde a_1$ -- Eq.~\eqref{subeq-1: coef of aeff nd_x}
		& $0$ & $-0.1756$ & $0.1192$ & $0.1756$ \\
		
		$\tilde a_2$ -- Eq.~\eqref{subeq-2: coef of aeff nd_x}
		& $0.4055$ & $0.2851$ & $0.3580$& $0.2851$  \\
		
		$\tilde a_3$ -- Eq.~\eqref{subeq-3: coef of aeff nd_x}
		& $0$ & $0.0658$ & $-0.0539$ & $-0.0658$ \\
		
		$\tilde a_4$ -- Eq.~\eqref{subeq-4: coef of aeff nd_x}
		& $0$ & $0.0036$ & $0.0020$ & $0.0036$ \\
		\bottomrule
	\end{tabular}
\end{table}

%

Thus, we can see that a fold, due to zero effective stiffness of the valve is what underlies the blowdown effect, that is designed into simple DSOPRVs. In turn,  it is the geometric shape of the valve through $A_{\rm eff}(x)$ (and to a less extent $A_{\rm ft}(x)$) that can give rise to 
these folds. 


To illustrate this effect, we have performed calculations for a basic 2J3 valve with parameters defined in \Cref{tab:parameters}.  
We then consider the variation of the 
effective area from a simple constant case $A_{\rm eff} =A_0$ to nonlinear cases of $A_{\rm eff}(x)$
by introducing different half-open angles and set pressures. 
\Cref{tab:-aeffs-from-different-open-angles} gives the calculated values of the dimensionless form of the coefficients of the quadratic form of 
$A_{\rm eff}(x)$ in Eq.~\eqref{eq:Aeff_Taylor} for four different values of
$\alpha = \pi -\varphi$. 
\begin{figure}
	\centering
	\hspace{-2em}
	\begin{tabular}{c c}
		\rotatebox{90}{\small Dimensionless lift, $x/(D/4)$}
		&
		\begin{tabular}{
				>{\centering\arraybackslash}m{6.2cm}
				>{\centering\arraybackslash}m{6.2cm}}
			
			\vspace{-2em}
			\begin{overpic}[width=5cm]{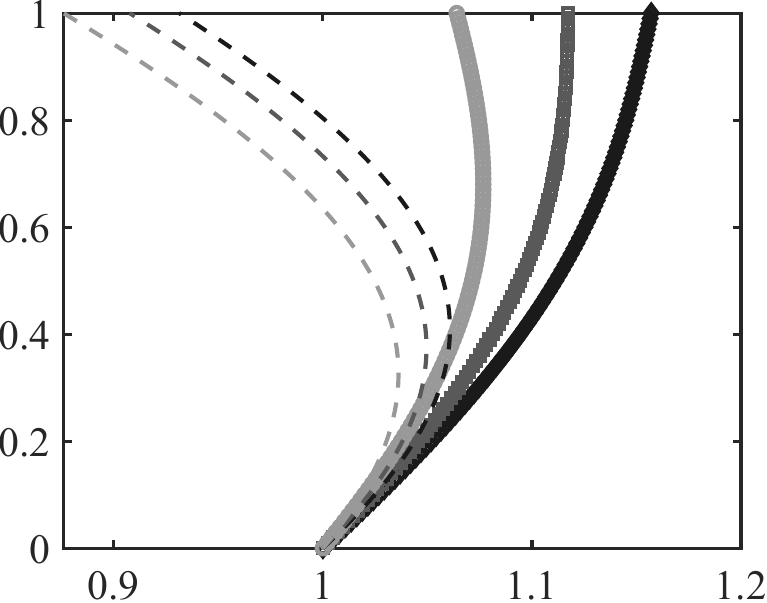}
				\put(-10,80){\normalsize{(a)}}
			\end{overpic}
			&
			\vspace{-2em}
			\begin{overpic}[width=5cm]{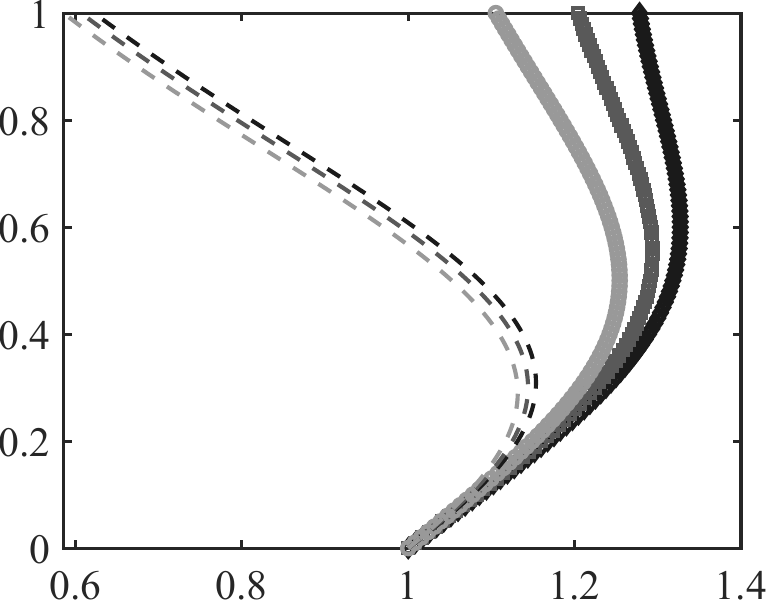}
				\put(-10,80){\normalsize{(b)}}
			\end{overpic}
			\\[1em]
			
			\begin{overpic}[width=5cm]{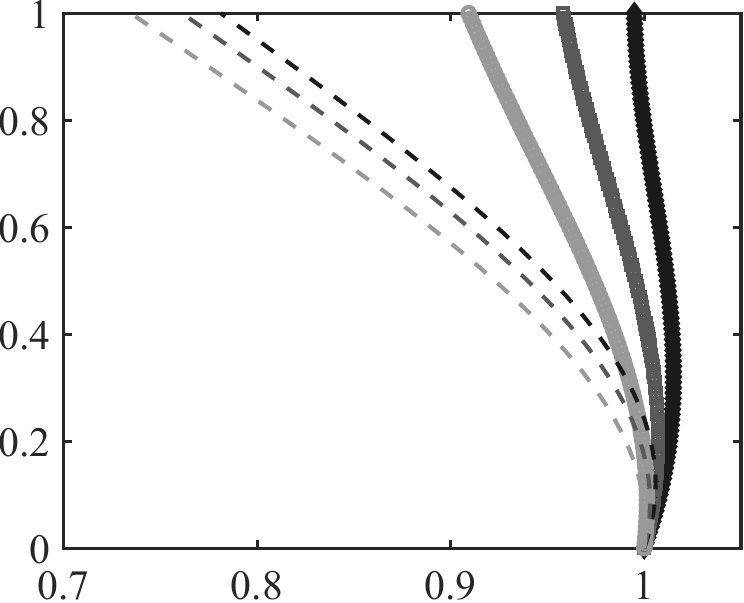}
				\put(-10,80){\normalsize{(c)}}
			\end{overpic}
			&
			\begin{overpic}[width=5cm]{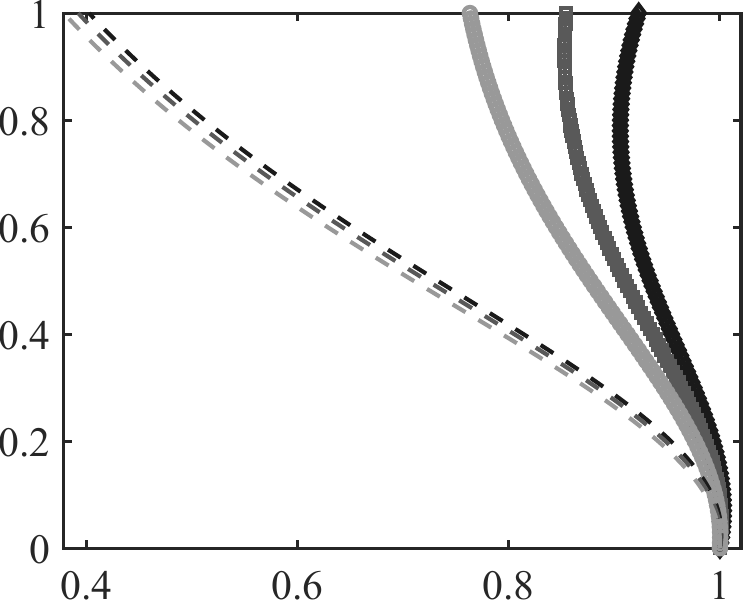}
				\put(-10,80){\normalsize{(d)}}
			\end{overpic}
			\\[1em]
			
			\begin{overpic}[width=5cm]{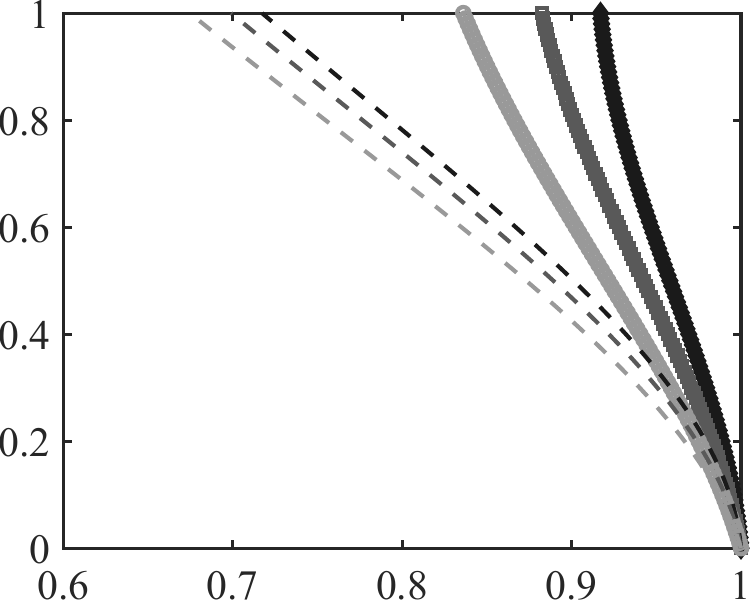}
				\put(-10,80){\normalsize{(e)}}
			\end{overpic}
			&
			\begin{overpic}[width=5cm]{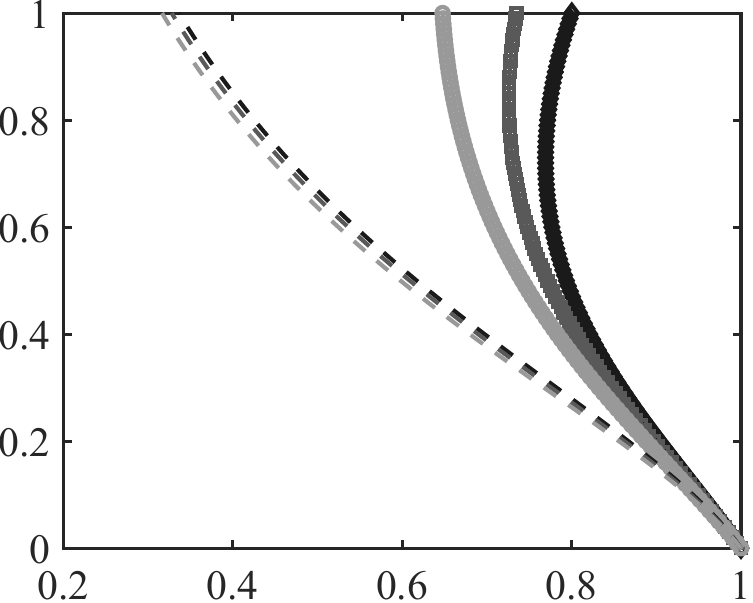}
				\put(-10,80){\normalsize{(f)}}
			\end{overpic}
			\\[1em]
			
			\begin{overpic}[width=5cm]{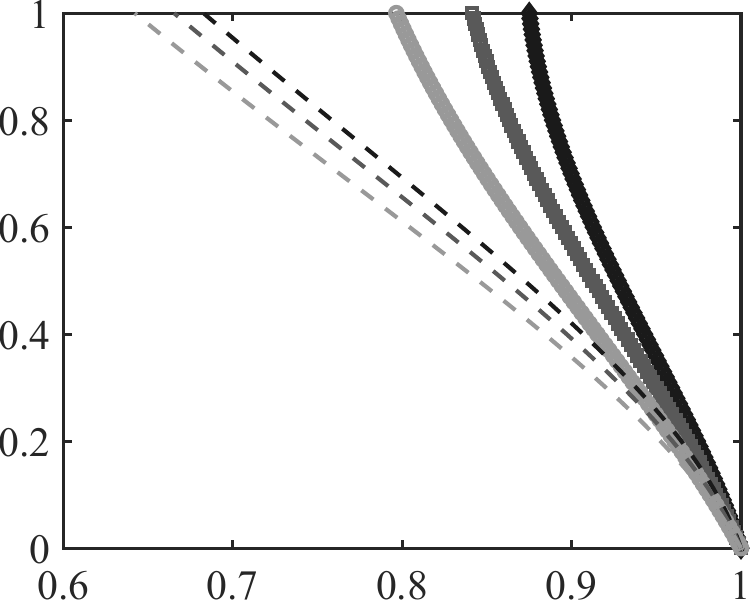}
				\put(-10,80){\normalsize{(g)}}
			\end{overpic}
			&
			\begin{overpic}[width=5cm]{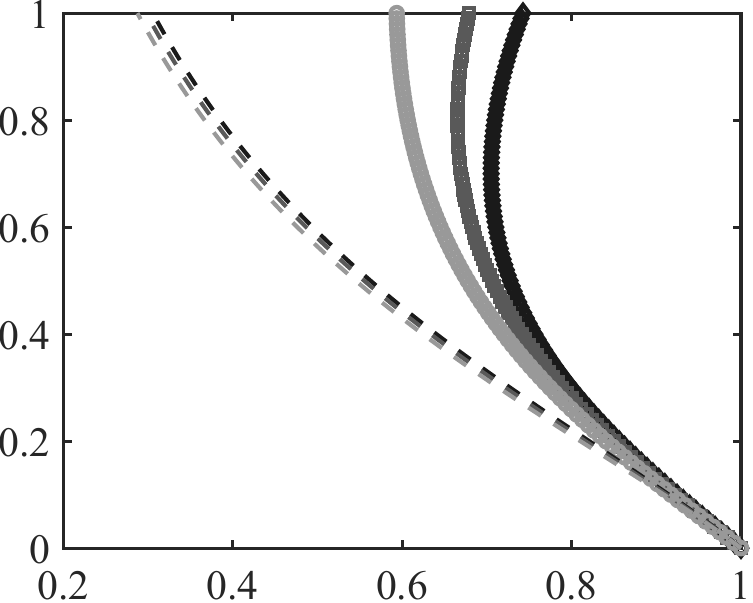}
				\put(-10,80){\normalsize{(h)}}
			\end{overpic}
			\\[-0.2em]
			\multicolumn{2}{c}{\small Dimensionless pressure}
		\end{tabular}
	\end{tabular}
	
	\caption{Equilibrium curves for different half-open angles $\alpha$ and service conditions, under varying set pressure. The left column corresponds to gas service and the right column to liquid service. The rows correspond to $\alpha=\pi/3$, $\pi/2$, $3\pi/5$, and $2\pi/3$, respectively.
		Solid and dashed lines show the dimensionless lift $x/(D/4)$ plotted against the tank pressure $(p_r-p_b)/p_{\rm set}$ and the valve-end pressure $(p_v-p_b)/p_{\rm set}$, respectively. The dimensionless parameter is defined as
		$\delta=p_{\rm set}/p_b$. Increasing line darkness denotes decreasing values of $\delta$, namely $\delta=15, 10, 5$.}
	\label{fig:geometric-aeff-demos}
\end{figure}

Figure \ref{fig:geometric-aeff-demos} shows the corresponding valve characteristics for each geometry,
in both liquid and gas service in each case. Fold points are not always obvious when lift is plotted against the valve pressure $p_v$, shown by dashed lines, but are more readily observed when lift is plotted against the reservoir pressure $p_r$, shown by solid lines.  

\begin{figure}
	\centering
	\subfloat[]{
		\includegraphics[width = 0.45 \textwidth]{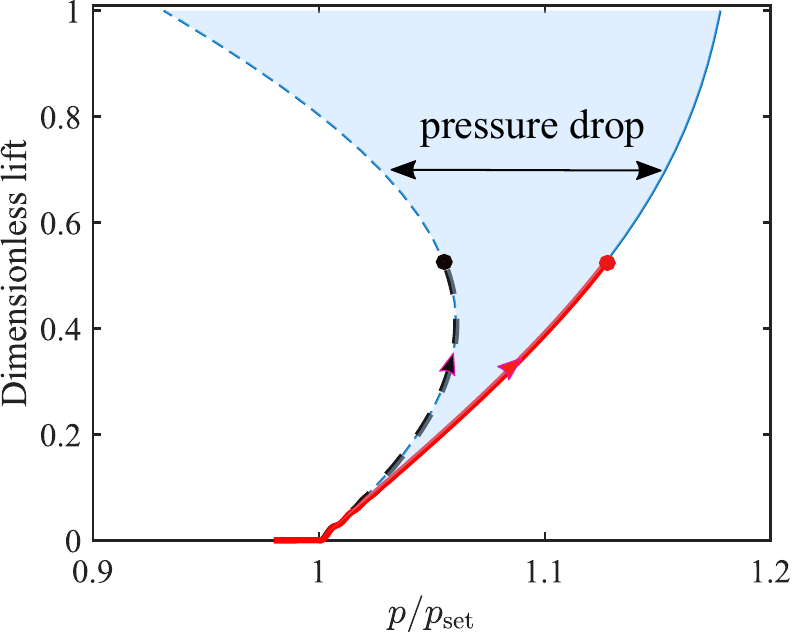}
	}
	\subfloat[]{
		\includegraphics[width = 0.45 \textwidth]{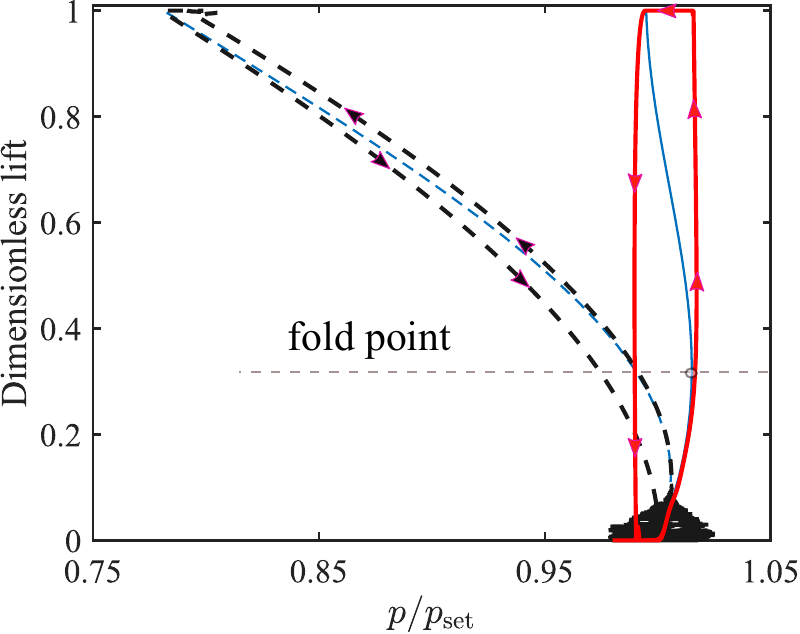}
	}
	\caption{Solutions to the PDM  superimposed onto the different equilibrium curves, for the cases in 
		Fig.~\ref{fig:geometric-aeff-demos} with: left, $\alpha = \pi/3$; right, $\alpha = \pi/2$. 
		Simulations were performed  with $\dot{m}_{\rm in} = 0.5~{\rm kg/s}$ and  $L = 0.1 ~{\rm m}$. Solid lines denote lift plotted against tank pressure, while dashed lines denote lift plotted against valve pressure. Thin lines denote equilibrium solutions and thick lines denote simulated trajectories. }
	\label{fig:blow-down-effect-show-case}
\end{figure}

Figure~\ref{fig:blow-down-effect-show-case} gives more details of the difference between a cone and a disk valve by using the PDM for a gas service valve with $\alpha=\pi/3$ and $\alpha=\pi/2$. In each case we provide
time history responses for a whole loop of valve opening and closing. 
To enable this to happen, we slowly increase the parameter $\dot{m}_{\rm in}$
up to a maximum value and then decrease it. The simulated results are superimposed on top of the calculated valve equilibrium curves in each case. We also plot the results in terms of both the tank pressure $p_r$ and the valve pressure $p_v$.

For the cone-shaped valve, we see from Fig.~\ref{fig:blow-down-effect-show-case}(a) that the simulated dynamics lies 
on top of the equilibrium curve. This is to be expected because 
$k_{\rm eff}$ (which we can interpret as the slope of the 
$x$-versus-$p_r$ curve) is strictly positive and hence the equilibrium is
stable for all lift-values. 

The dynamics for the disk valve in Fig.~\ref{fig:blow-down-effect-show-case}(b) is more complex. Here we see that after opening, the simulation hugs the equilibrium curve only up to the fold point in the $x$-versus-$p_r$ curve. There is then a rapid jump from this fold point up to full lift. 
As the pressure is then decreased, there is a pseudo equilibrium, which persists  down to a value of $p_r$ which is a little bit lower than 
$p_{\rm set}$ (by about 1\%).  
At this point there is a boundary equilibrium
because the unstable portion of the equilibrium curve reaches full lift.
At this point there is another jump, this time down to valve closure. 

We note that this jump behavior is much more easily understood by looking at 
$p_r$ (the black trajectories and red characteristic curve) than at $p_v$ (the green trajectories and blue characteristic curve). 
However, there are several interesting features that can be observed  when looking at the lift versus valve pressure. 
First we note that at large lift, there is a considerable inlet pressure loss
(the difference between $p_v$ and $p_r$). Second, note how the valve pressure time series show a significant amount of {\em chattering}, when the valve jumps into closure, that is a geometrically decreasing sequence of impacts in which there is a small amount of flow through the valve, leading to fluctuating valve pressure. Note though that the reservoir pressure seems largely unaffected by these oscillations.

In summary, these calculations demonstrate that the blowdown effect can indeed be explained by the presence of fold points leading to negative stiffness on the equilibrium branch of lift $x$ versus tank pressure $p_r$.
In turn, these fold points arise from the net effect of 
$A_{\rm eff}(x)$ and effect of inlet pressure drop. 
We note that the case $\alpha=\pi/2$ 
demonstrates only a small blowdown (about 1\%). As we can see from the shape
of the equilibrium curves in Fig.~\ref{fig:geometric-aeff-demos} though, 
blowdown is much more pronounced for larger $\alpha$-values.

%

\subsection{Quarter-wave dynamic instability}
\label{sec:reduced-order-model}

As reviewed in the Introduction, it is well established that 
statically stable dynamic equilibrium operations of  DSOPRVs can be subject to a dynamic flutter instability. Here 
the dynamics of the valve act as a negative damper to oscillatory acoustic modes of  velocity and pressure perturbation to the mean flow in the inlet pipe. 
To that end an expansion is considered of the form 
\begin{align}
	\label{eq:trial-solution-coll-dimensional}
	\begin{split}
		p(\xi,t)  &= p(0,t)+ \sum_{j=1}^{N} B_j(t) \sin \left( \frac{\pi\xi(2j-1)}{2} \right), \\
		v(\xi,t)&= v(0, t) + \sum_{j=1}^{N} C_j(t) \cos
		\left( \frac{\pi\xi (2j-1)}{2} \right),
	\end{split}
\end{align}
which naturally satisfies the boundary conditions at both ends of the pipe. 
However, it has also been established 
\cite{HosCh2014-I,HosCh2015-II,HosCh2016-III,HosCh2017-IV,Ma_experiments,Schmidt_models} 
that, provided the pipe is not excessively long, the initial instability upon increasing the pipe length $L$ (or equivalently, reducing the mean pipe velocity) is the mode for $j=1$, which represents the fundamental quarter-wave mode.  

Thus, to understand
the initial dynamic instability, it suffices to truncate Eq.~\eqref{eq:trial-solution-coll-dimensional} after one term. This leads to the so-called quarter-wave model (QWM).
The derivation details are given in  \ref{sec:Appendix_III}, 
where we improve the derivation presented in \cite{Hos2014,BazsoCh2014,HosCh2015-II}, to include a general form of $A_{\rm eff}(x)$ and correct treatment of the inlet
pressure loss through the function $\chi$).
Setting $N =1$ in Eq.~\eqref{eq:trial-solution-coll-dimensional}  
\begin{subequations}
	\label{eq:dimensional-QWM-OG}
	\begin{align}
		m \ddot{x} & = (p_0 + B - p_b) {A}_{\rm eff}(x) - c \dot{x} - k(x + x_0),\\
		\dot{p}_r & = \frac{a^2}{V} \left [ \dot{m}_{\rm in} - \frac{\pi D^2}{4} \rho (v_L + C) \right ],\\
		\dot{B} &= - \sqrt{2} \dot{p}_0 + a^{2} \rho  \frac { \pi } { 2 L }C - (\sqrt{2} v_{L} + C) \frac
		{\pi}{2\sqrt{2}L} B,\\
		\dot{C} &=  -\sqrt {2}\dot{v}_L + ( \sqrt { 2 } 
		v_{ L } + C ) \frac { \pi } { 2 \sqrt{2}L } C -\frac { 1 } {
			\rho } \frac { \pi } { 2 L } B  { -}\lambda \frac { L } { D } ( \sqrt { 2 } v_L+ C ) \left| \sqrt {
			2 } v_L + C \right |.
	\end{align}
\end{subequations}

\begin{figure}[t]
	\centering
	\includegraphics[width = 0.5 \linewidth]{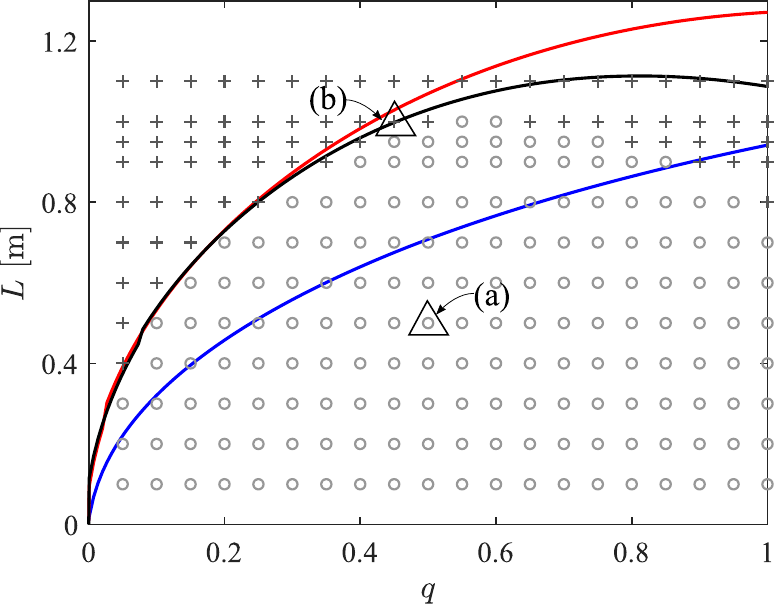}
	\caption{Stability chart showing the onset of 
		valve flutter via a quarter-wave instability for the 2J3 valve in gas service for different pipe lengths and relative  mass flow rates $q = \dot{m}_{\rm in}/ \dot{m}{\rm in,cap}$.
		The markers represent the results of the full PDM simulation of valve opening: crosses denote unstable cases and circles stable cases. The labels (a) and (b) refer to the two cases where the simulations are presented in detail in Fig.~\ref{fig: PSTool example-q0.5}. The 
		three solid lines represent the 
		Hopf bifurcation boundary computed using the QWM. The blue line is the original analytical prediction from \cite{HosCh2015-II}. The red line is the corresponding computed curve that includes both the nonlinear convection terms but not inlet pressure loss. The black line includes the inlet pressure drop $\chi$.}
	\label{fig: comparison between CO and ANA}
\end{figure}

By ignoring the inlet pressure drop effect, setting $\chi=0$, the 
algebraic complexity of the model is hugely reduced because $p_0 = p_r$. In this case, an approximate analytical formula for the quarter-wave instability threshold was derived in \cite{HosCh2015-II,HosCh2017-IV}. 
This precise analytic condition is limited however, because it was based on neglecting the pipe inlet pressure drop, assuming constant
$A_{\rm eff} = A_0$, linear $A_{\rm ft}$, ignoring the (nonlinear) convection terms in the model and neglecting pipe friction. 
Of these simplifications, the results in Fig.~\ref{fig:blow-down-effect-show-case} suggest that ignoring the inlet pressure drop can have the greatest effect. 
Nevertheless, it remains relatively straightforward to find the quarter-wave instability boundary for the full version of the QWM Eq.~\eqref{eq:dimensional-QWM-OG} that includes all of these effects. (Note though that in our computational results in both the PDM and QWM we have assumed the pipe friction coefficient $\lambda=0$). 

To do this, one simply needs to compute the 
dynamic equilibrium $x=x_e$, $p=P_r(x_e)$, $B=C=0$. Then the  
Eq.~\eqref{eq:dimensional-QWM-OG} can be linearized about this point and the appropriate
Jacobian matrix formed. See \ref{sec:Appendix_III} for the details. The condition for linear stability is then the usual one that all the eigenvalues of this matrix should be in the left half-plane. Flutter instability arises upon varying a parameter (such as the constant mass inflow $\dot{m}_{\rm in}$) when a pair of pure imaginary eigenvalues cross into the right half-plane. In the dynamical systems literature, such an 
instability is known as a {\em Hopf bifurcation}. The process of tracking the locus of such Hopf bifurcations can be automated using 
numerical continuation software (see e.g.~\cite{HosCh12gr,BazsoCh2014} for examples of using such numerical bifurcation analysis on reduced-order valve models). 

\begin{figure}
	\centering
	\subfloat[$q = 0.5$, $L = 0.5$m.]{
		\centering
		\includegraphics[width = 0.49 \textwidth]{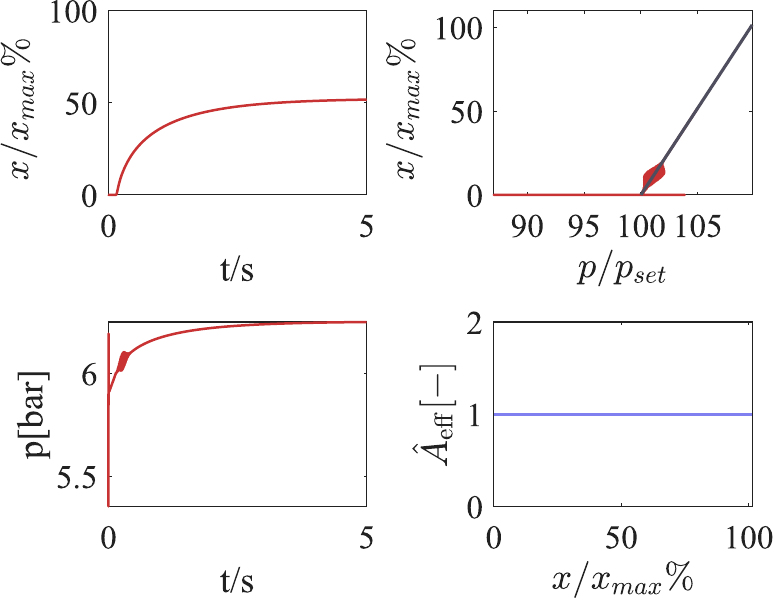}
	}
	\subfloat[$q = 0.5$, $L = 1$m.]{
		\centering
		\includegraphics[width = 0.49 \textwidth]{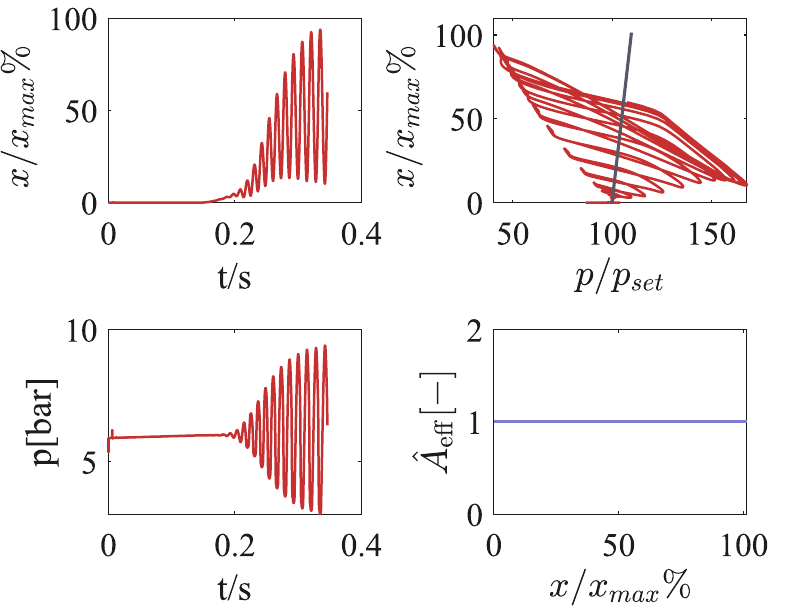}
	}
	\caption{Example runs of PDM at the two labeled points in Fig.~\ref{fig: comparison between CO and ANA}. In this and subsequent similar figures, red lines show the simulation trajectories, and black solid lines are analytical approximations to the valve equilibrium lift $x_e$ versus tank pressure $p_r$. In each case, the upper left panel shows the time history of valve lift as a fraction of $x_{\rm max}$;  the upper right panel shows the lift versus tank pressure $p_r$; the lower left shows valve pressure $p_v$ against time; the bottom right depicts 
		${A}_{\rm eff}  = A_0 \hat{A}_{\rm eff}$ 
		which for this case is trivial because $\hat{A}_{\rm eff}  = 1$. }
	\label{fig: PSTool example-q0.5}
\end{figure}

We have carried out such a process for the same example 2J3 valve in gas service used in \cite{Hos2014}, under variation of the mass flow rate $\dot{m}_{\rm in}$ and pipe length $L$,
and compared the results with full simulations using the PDM. The parameters values used are given in \Cref{tab:parameters}. The results are presented in Fig.~\ref{fig: comparison between CO and ANA}.
For each full simulation we start with the valve closed and the $p_r$ below set pressure. Then we allow the pressure to increase and the valve to open. Examples of such simulations are given in 
Fig.~\ref{fig: PSTool example-q0.5}. 

Looking at the results in 
Fig.~\ref{fig: comparison between CO and ANA} we see that all of the curves computed using the QWM follow the same trend as the full simulation results. That is, stability is achieved for sufficiently short pipe lengths $L$ which varies approximately quadratically with mass flow rate, so that instability occurs for
long pipes and low flow rates. However, we see that for long pipes, inclusion of pipe friction and convection (the red curve in 
Fig.~\ref{fig: comparison between CO and ANA}) is necessary to get the correct quadratic approximation to the true instability boundary. This curve 
provides a much better approximation to the true stability boundary for low relative mass flow rates $q\leq 0.4$ than the original analytical approximation (blue curve). However, if we do not include the inlet pressure drop, then the QWM greatly over-predicts the pipe length corresponding to instability for flow rates that are closer to capacity. Here we see that inclusion of inlet pressure drop $\chi$ (black curve) captures the  true trend for the instability curve for higher $q$.

It is worth noting though that there is still a small discrepancy between the QWM prediction and simulation results, which becomes more extreme close to capacity 
flow rates. We should note here though that the stability curve and the simulation results are measuring different things. Note that the stability curve derived from the full QWM is actually the true curve of quarter-wave-instability for the full PDM because it is derived by seeking Hopf bifurcations 
corresponding to infinitesimal perturbations in the direction of the exact quarter-wave mode of the pipe. 
However, the simulation is for transient results upon opening the valve. For higher $q$ these transients 
represent a large perturbation from the equilibrium
valve lift. In the parameter regime where there is discrepancy between the transient simulations and the 
linear stability results, we find that the perturbation associated with valve opening is too large
to be attracted to the stable equilibrium lift. Instead, the system finds another attractor, typically one involving non-decaying chattering. 

In fact, even in the stable case depicted at point (a) 
in Fig.~\ref{fig: comparison between CO and ANA}, which is well inside the stability region, we see from the simulation results in Fig.~\ref{fig: PSTool example-q0.5}(a) that there are small oscillations upon valve opening. However these oscillations quickly decay. For the unstable case
Fig.~\ref{fig: PSTool example-q0.5}(b) we see instead that the opening transients do not decay. Instead they grow and shortly after $t=0.3$ become so extreme
that the valve pressure becomes 0 at which point the flow through the pipe would reverse.


\section{Analysis of enhanced stability valves}
\label{sec:stability}

The purpose of this section is to analyse  of the enhanced stability concept introduced in Fig.~\ref{fig:stopper_below_fold} of \cref{sec:design-concept-intro}. The key idea there is to design the blowdown curve and upper stop location to avoid the Hopf bifurcation associated with the quarter-wave instability. Here, there would be no stable dynamic equilibrium that is accessible to the valve. Instead, it would operate at the pseudo equilibrium at the upper stop. Note that the corresponding dynamic equilibrium at these lift values would, in theory, be susceptible
(violently) unstable to flutter instability, so we need to rely on the
pseudo equilibrium being sufficiently stable, and the transient opening to be sufficiently rapid to get to the upper stop before instability sets in. 

The case illustrated in Fig.~\ref{fig:stopper_below_fold}, is for a 2J3 valve with nontrivial $A_{\rm eff}(y) = 1 + y^2$ corresponding to 
$x_{\rm stop} =0.4 x_{\rm max}$ and $e_1=0.2$. The thin red lines in the figure show the simulated trajectory using the PDM.
Here we do indeed see that the
pseudo equilibrium is quickly reached after the valve is opened and instability is avoided. We now need to demonstrate this effect more generally and determine how the extent of the enhanced stability (compared with an unmodified valve) depends upon system parameters. All computations are presented for modifications to a 2J3 valve in gas service. 

From a practical point of view, we can think of two equivalent ways of achieving this enhanced stability:
\begin{enumerate}
	\item Strategy I -- from a system designer's point of view --- choose a larger capacity valve and restrict its lift through the location of $x= x_{\rm stop}< x_{\rm max}$, as in Fig.~\ref{fig:stopper_below_fold};
	\item Strategy II -- from a valve manufacturer's point of view ---
	while keeping $x_{\rm stop} = x_{\rm \max}$, manipulate the valve's geometry so that the lift-versus pressure characteristic has its fold point (at the red diamond in
	Fig.~\ref{fig:stopper_below_fold}) for $x>x_{\rm max}$.
\end{enumerate}
From a mathematical analysis point of view, these two approaches are equivalent, because in each case all that matters is  that equilibrium lift versus tank pressure curve should have  negative slope all the way up to the upper stop. 
Therefore, in what follows we do not distinguish between $x_{\rm stop}$ and
$x_{\rm max}$. 

\begin{figure}
	\centering
	\vspace{-0.5cm}
	\subfloat[$q = 0.5; \; L = 2 m.$]{
		\centering
		\includegraphics[width = 0.49 \textwidth]{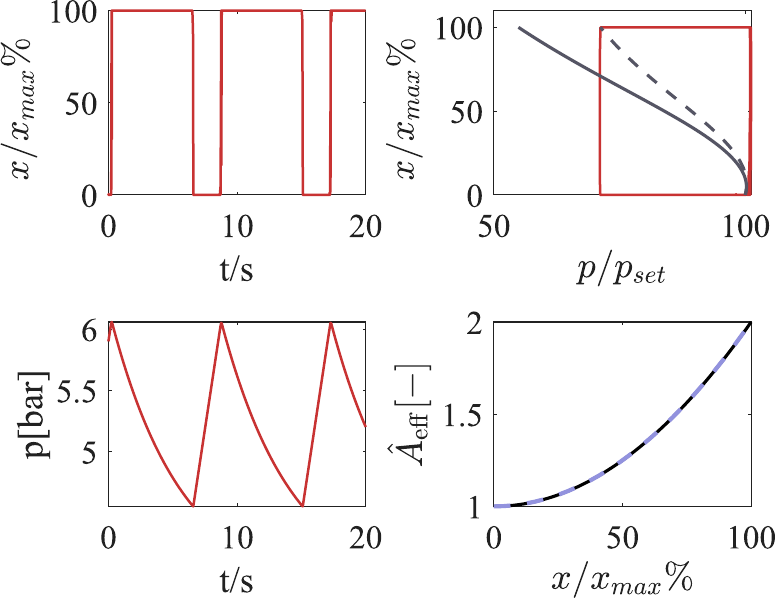}
	}
	\subfloat[$q = 0.5; \; L = 5 m.$]{
		\centering
		\includegraphics[width = 0.49 \textwidth]{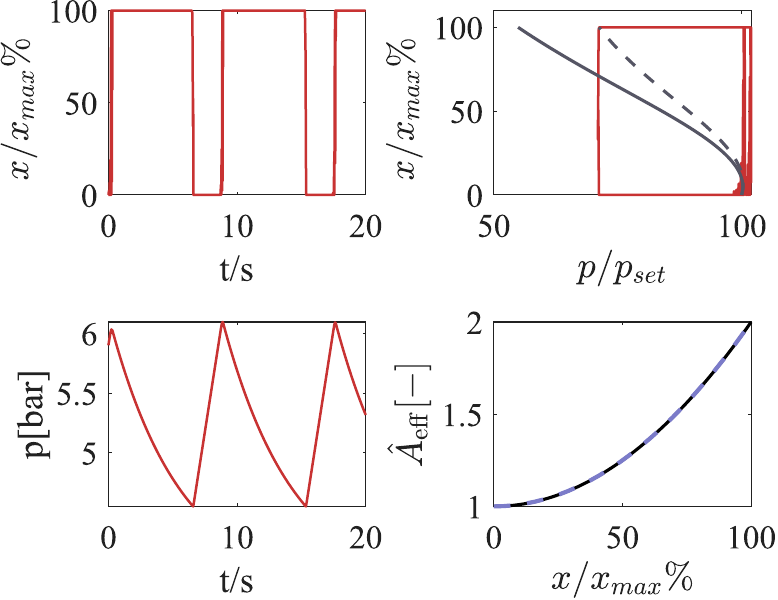}
	}
	\caption{ Similar to 
		Fig.~\ref{fig: PSTool example-q0.5} but for $\hat{A}_{\rm eff} = 1 + y^2,$  and $r =0.2$. Here, in addition, the black solid/dashed line in the upper right plots is the (unstable) equilibrium lift versus pressure 
		$p_r$/$p_v$.}
	\label{fig: PSTool examples - blowdown}
\end{figure}

\begin{figure}
	\centering
	\hspace{-2.5em}
	\subfloat[]{
		\includegraphics[width = 0.495 \textwidth]{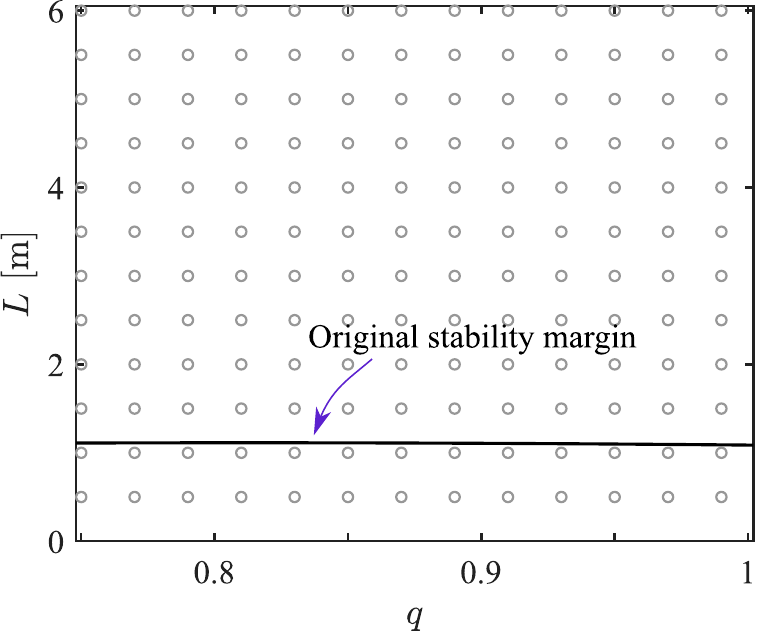}
		\label{fig: stability_enhanced_diagram}
	}
	\subfloat[]{
		\includegraphics[width = 0.5 \textwidth]{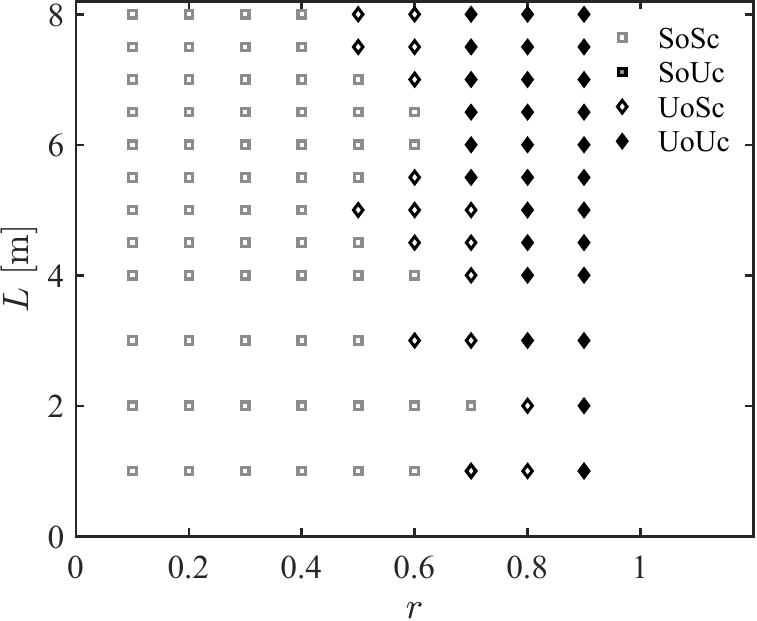}
		\label{fig:stability_chart_q0.5}
	}
	\caption{(a) Stability diagram obtained by running the PDM with the same parameters as in \cref{fig: comparison between CO and ANA} but with  $\hat{A}_{\rm eff} = 1 + y^2,$ for varying $L$ and $q$ with coefficient of restitution $e_{0,1}=r=0.2$. Superimposed on this is the QWM  (b) The stability chart when the $q=0.5$.  Here 
		the key shows cases corresponding to: stable opening and stable closing (SoSc),  stable opening and unstable closing (SoUc),  unstable opening and stable closing (UoSc), and  full chattering (UoUc) }
\end{figure}

To investigate the validity of our approach more generally, we consider a case of a disc with collar angled at $\alpha = 2\pi/3$ (see row (h) in Fig.~\ref{fig:geometric-aeff-demos} )
which leads to a quadratic
$A_{\rm eff}(x)$ curve. The results are shown in Fig.~\ref{fig: PSTool examples - blowdown}
and \ref{fig: stability_enhanced_diagram}. These results can
can be compared with the corresponding ones for the  
original 2J3 valve ((which had $\hat{A}_{\rm eff}=1$ and $\alpha= \varphi = \pi/2$)
in Fig.~\ref{fig: PSTool example-q0.5} and
Fig.~\ref{fig: comparison between CO and ANA}.

\Cref{fig: PSTool examples - blowdown} shows simulation results for 
mass inflow rate is s$\dot{m}_{\rm in} = q \dot{m}_{\rm in,cap}$ with $q=0.5$. and two different inlet pipe lengths $L$ that would be beyond the stability limit for the unmodified
valve. The time simulations show that the valve opens to full lift in both cases once the tank set pressure $p_{\rm set}$ is reached. Because the inflow rate is less than capacity, there is no further increase in pressure at this point. Instead, the pressure
slowly decreases and the pseudo equilibrium branch is tracked until we reach the point
at which blowdown occurs, which is at $p_r=P(x_{\rm stop})$, where $P(x_e)$ is the dynamic  equilibrium characteristic defined in Eq.~\eqref{eq:static_curve_pr_x_liquid}.
For this example and flow rate, we can  calculate from the shape of the
corresponding $A_{\rm eff}(x)$ that the blowdown is $\Delta_{\rm bd}$ is 28\%. Here we define blowdown via
\begin{equation}
	\label{eq:blowdown}
	\Delta_{bd} = \frac{P(x_{stop}) - p_{\rm set}}{p_{\rm set}}.
\end{equation}

At this blowdown pressure we have a so-called boundary equilibrium bifurcation,
which takes the form or a {\em nonsmooth fold} \cite{DiBernardo2008,TangSimpson}. 
For lower pressures, the pseudo equilibrium no longer exists as 
$\ddot{x_{\rm stop}}$ is now  negative. As $p_r < p_{\rm set}$, there is no stable dynamic equilibrium either so the valve rapid shuts. The dynamics now then proceed with the
the valve closed and the pressure rapidly rising, until the set pressure is once again reached and this cycle begins again. Note this valve cycling between open and closed as the pressure is not at all rapid --- the is frequency $\sim 0.1$ Hz, whereas the damaging chattering due to quarter-wave instability typically has frequency $\sim 100 \rm{Hz}$.

The simulations show there is no visible evidence of
transient chattering of the valve either on openning or closing. The
key here is that we used a low enough coefficient of restitution of
$e_0 = e_1 = 0.2$.
Hence each impact with the upper stop dissipates energy. This enables the
dynamics to be quickly attracted to the pseudo equilibrium.

We have repeated the simulations in Fig.~\ref{fig: PSTool examples - blowdown} for a large range of flow rates and inlet pipe lengths. 
The results are shown in 
Fig.~\ref{fig: stability_enhanced_diagram},
where we have zoomed in on the higher mass flow rates.  We find qualitatively the same result as in
Fig.~\ref{fig: PSTool examples - blowdown}
for all flow rates $q$ up to capacity, and for all pipe lengths $L$ that we have simulated. For comparison, the solid line in the figure shows the stability threshold (maximum pipe length $L$ for the unmodified valve.

To investigate the influence of coefficient of restitution on the stability of our
modified valve, we have repeated the computation in Fig.~\ref{fig: stability_enhanced_diagram}, but now keeping the flow rate fixed and allowing the coefficient of restitution $r$ to vary, along with the pipe length $L$. Now we find that the stability area is
greatly reduced as we take values of $r=1$ (which corresponds to no energy dissipation at impact). As the diagram shows there are three different kinds of instability that we observe, depending on whether non-decaying chatter is triggered on valve opening, valve closing or both.

\section{Conclusion}\label{sec:discussion}

This paper has proposed and analyzed a key new design concept for
enhancing stability of DSOPRVs  against oscillatory
flutter or chatter instabilities. Namely, partial valve opening is a
bad thing and should be avoided. As the valve opens to full lift, then
a stable equilibrium can be reached by providing an energy dissipation
mechanism, which we have represented via a simple coefficient of
restitution.

To arrive at this concept, we have produced a novel explanation of the
origins of the blowdown effect in terms of folds in the shape of the
equilibrium valve lift-versus-pressure characteristic.  This has
enhanced our earlier work, summarized in \cite{HosCh2017-IV} by
showing explicitly how the effective-area-versus-lift characteristic
can be used to define an effective stiffness, the vanishing of which
leads to the fold. The key idea to ensure that full lift is reached
upon valve opening, is that this fold curve should occur for
unrealistically high valve lifts that are {\em above} the upper
stopper $x_{\rm stop}$. Moreover, if we have the ability to modify the
valve pressure-versus lift characteristic $P(x_e)$ as in 
Eq.~\eqref{eq:static_curve_pr_x_liquid}, 
then a specified amount of blowdown can be designed through the formula
Eq.~\eqref{eq:blowdown}.

While this paper is not primarily aimed at practical design advice, we
have proposed two different strategies to enhance the stability of
standard, so called \emph{proportional-lift}, DSOPRVs. Both essentially
involve redesign of the lift-versus-pressure characteristics of the
valve through modification of the effective area versus lift function
$A_{\rm eff }(x)$. These we can think of as retrofit (Strategy I) or
redesign (Strategy II). In the first case, the engineer would select
and oversized valve for the job at hand (i.e.~one that would have too
large a capacity) and restricting its lift.  There is recent experimental evidence that this strategy is likely to be beneficial through the work of Smith \& Desai \cite{Smith2025}
who performed hundreds of laboratory tests on valves with restricted lift. They found that stability was greatly enhanced. 

The second strategy is
aimed at valve manufacturers to redesign the internal geometry of
their valve to ensure an optimal shape of the effective area versus
lift curve. Our CFD results show that this latter approach may also be
thought as optimization of the discharge jet angle $\phi$. However,
jet angle alone is a rather crude way of characterizing the complete
shape of $A_{\rm eff}(x)$.

We should also stress other ways this work also improves the state of
knowledge about DSOPRVs generally. First, we have produced an
improvement to our previously published model of a pressure-relief
system, incorporating a valve and its input pipeline, both for the
case of gas and liquid service. We believe we have now presented this
model in the clearest manner yet, which should be of benefit to future
researchers.  Second, we have improved our previously-published
quarter-wave reduced-order model. By correctly including inlet
pressure loss and pipe friction, we have shown that this model is able
to predict the instability threshold (of unrestricted valves) much
more accurately. Finally, we have used CFD to provide validation of
our use of effective area versus lift curves as a good way to capture
a valve's characteristic.

There are many limitations to the work presented. First and foremost,
while our PDM model has been fully calibrated with experimental
measurements in \cite{HosCh2014-I,HosCh2016-III}, we have yet to
directly demonstrate the effectiveness of our new enhanced stability strategy in laboratory experiments. We note however
the aforementioned work \cite{Smith2025} that is highly promising, but did not explicitly consider each valve's blowdown characteristics. This paper is already long, and direct
experimental validation is left for future work.

Another limitation is that although the concept is quite general, we have
only done detailed computations for a valve that is servicing an ideal
gas. We have performed some preliminary computations for the liquid
service case, which are promising, but have found that more care is
needed because of the more violent {\em subcritical} version of the
quarter-wave Hopf bifurcation in this case
\cite{HosCh2016-III,HosCh2017-IV}. A detailed examination of this case
will be the subject of future work. Future work should also consider a more
detailed global stability analysis using ideas from non-smooth
dynamical systems theory, including the notion of whether the valve
opening transient is within the basin of attraction of the
pseudo-equilibrium.

Another area ripe for further investigation concerns the role of reservoir pressure. Our results show that the tank-pressure-versus-valve-lift characteristic, $p_r=P(x_e)$, is the key quantity for understanding valve stability. We have seen though that there can be a large difference between $p_r$ and the pressure $p_v$ at the valve end of the pipe for high flow rates.  Note that $p_v$ versus lift is more likely to be the known characteristic of any valve, whereas determination of $p_r$ requires estimation of losses in the inlet pipe.
In particular, 
pipe inlet pressure loss $\chi$ plays a vital role, and our results
(see e.g.~Fig.~\ref{fig:geometric-aeff-demos} and  \ref{fig: PSTool examples - blowdown})
show that ignoring it leads to false conclusions
about blow-down. Blowdown would typically occur for much lower pressures if the blowdown formula Eq.~\eqref{eq:blowdown} were designed using $p_v$ rather than $p_r$. Also, our simulations have so far been run without pipe friction loss, that is
$\lambda=0$. For relatively short inlet pipes $L \sim$ 1--5 m we have found in earlier work that these frictional losses are typically much smaller than the inlet pressure loss. However, given that our enhanced stability valves show stability for longer inlet pipes, for sufficiently small $r$, the effect of friction also needs to taken into account
when calculating $p_r$ from $p_v$.

We have also not considered more complex fluids such as non-ideal gas, steam,
non-Newtonian fluids, multiphase flow etc. Nor have we considered practical
implementation details. For example, given the importance we have
shown of sufficient damping associated with impact at the upper stop,
a simple restitution law is probably too simplistic. Practical
implementation should consider how to achieve such energy dissipation
via an additional vibration absorber. Finally, reported problems
associated with valve flutter and chatter are not restricted to the
enhanced stability concept we have proposed is likely to apply to
other valve topologies, beyond the DSOPRV,

In summary, while there remains much work to be done, we hope the
enhanced stability concept introduced in this paper will provide a
path for practitioners to finally eliminate damaging valve chatter, a
problem that has plagued industry since the design of the
pressure-relief valve.

\subsection*{Acknowledgments} 
The authors thank Ken Paul, Mike McNeely, Katherine Si, John Burgess and Dustin Smith for useful conversations that helped shape this work. Hong Tang was supported by the University of Bristol \& Chinese Scholarship Council joint studentship, No.202006120007.

%% file: appendix.tex
\section{Analytical approximation to valve discharge characteristics}\label{sec:Appendix_I}

This appendix provides more details on the improved derivation of valve discharge 
characteristics from that in 
\cite{HosCh2016-III} and \cite{HosCh2014-I} 
for the case of liquid and gas flows, respectively

\subsection{Flow rates through the valve}
We first summarize the derivation of formulae for the  output mass flow rate $\dot{m}_v$
which leads to the expression for the inflow $v_L$ through the valve.

For liquid flows we can calculate the change in total pressure to be 
$$
p_v = p_b + \frac{\rho}{2} v_{\rm ideal }^2, 
$$
where $v_{\rm ideal}$ is the theoretical velocity used to apply Bernoulli's equation. This theoretical velocity  neglects turbulence, viscosity and internal friction --- which effects are captured elsewhere in the calculation through the discharge coefficient. 
Bernoulli's principle then gives 
$$
v_{\rm ideal} = \sqrt{\frac{2}{\rho}(p_v - p_b)} 
$$
Then, using the usual definition of discharge coefficient, we get 
\begin{align}
\dot{m}^{\rm (liquid)}_{\rm v} & = C_d \rho A_{\rm ft}(x) v_{\rm ideal} \nonumber \\
  & = \sqrt{2} C_d A_{\rm ft} \sqrt{\frac{p_v - p_b}{\rho}}
\label{eq:m_v_liquid}
\end{align}

In the case of gas service valves, the flow through a narrow valve would typically be choked. Thus the output mass flow rate depends only on downstream pressure and we can write an 
expression for the outgoing mass flow rate using the usual 
formulae for gas dynamics through a nozzle with area
$A_{\rm ft}$ namely 
\begin{equation}
\dot{m}^{(\rm gas)}_v = C_d A_{\rm ft} \sqrt{\rho p_v \kappa
\left ( \frac{2}{1+\kappa} \right)^{\frac{\kappa+1}{\kappa-1}}
} 
\end{equation}
Here $\kappa$ is the usual adiabatic gas constant $\kappa$ (sometimes called $\gamma$) defined is the ratio of heat capacities at constant pressure and volume. 

In either case we calculate the 
magnitude of the 
corresponding discharge flow velocity as
\begin{equation}
v_{\rm out} = \frac{\dot{m}_v}{\rho A_{\rm ft} },
\label{eq:v_out} 
\end{equation}
which we assume exits at the jet angle $\varphi$.

Finally to calculate the formula Eq.~\eqref{eq:BC-V-at-1} 
for the inflow velocity 
$v_L$ to the valve we substitute the expression 
Eq.~\eqref{eq:v_out} into a simple mass balance equation: 
$$
\rho A_0 v_L = \rho A_{\rm ft} v_{\rm out} .
$$

\subsection{Effective area versus lift}

We next proceed to obtain an analytic expression for the effective area
characteristic. From Eqs.~\eqref{eq:BC-V-at-1} and 
\eqref{eq:Ffluid} we have 
$$
A_{\rm eff}(x) := A_0 +  \frac{1}{p_v-p_b} \dot{m}_v(v_{L} + v_{\rm out} \cos\varphi).
$$
Thus, we can substitute from the above expressions for $v_L$ and $v_{\rm out}$
to obtain 
$$
A_{\rm eff}(x) = A_0 +  C_d^2 C_\kappa^2 A_{\rm ft}(x) \frac{p^*}{p_v-p_b} \left(  \frac{A_{\rm ft}}{A_0} + \cos(\varphi) \right) 
$$
where we recall that $p^*=p-p_b$ for liquids and $p^*=p_v$ in the case of gas flow. 

Note that $p_b\ll p_v$ so, even in the case of gas flow we can write to high order of accuracy that 
$$
A_{\rm eff}(x) = A_0 \hat{A}_{\rm eff}(x),  
$$
where 
\begin{equation}
	\hat{A}_{\rm eff}\coloneqq  \left(1 + C_d^2 C_{\kappa}^2 \frac{ A_{\rm ft}(x)^2}{A_0^2} (1 +
	\frac{A_0}{A_{\rm
	ft}(x)} \cos\varphi)
	\right)
	\label{eq:hat_Aeff}
\end{equation}
which is exact in the case of liquid flows. 
Then we can substitute
the quadratic-in $x$ expression for the flow-through area 
Eq.~\eqref{eq:cone-shape-Aft}, to obtain an explicit analytic expression for $\hat{A}_{\rm eff} (x)$:
\begin{align}
	\hat{A}_{\rm eff} &= 1 +  a_1 x +   a_2 x^2 + a_3 x^3 +  a_4 x^4  
\end{align}
and in dimensionless valve lift $y = 4x/D$
\begin{align}
	\hat{A}_{\rm eff} &= 1 +  \tilde a_1 y +   \tilde a_2 y^2 + \tilde a_3 y^3 +  \tilde a_4 y^4  
\end{align}
where
\begin{subequations}
	\label{eq:hat_Aeff_coeffs}
	\begin{align}
		\tilde  a_1  &=\frac{ C_\kappa^2 C_d^2 \sin2\varphi}{2}, \label{subeq-1: coef of aeff nd_x}\\
		\tilde  a_2  &= \frac{C_\kappa^2 C_d^2 \sin^2\varphi }{4}(4 - \cos^2\varphi), \label{subeq-2: coef of aeff nd_x} \\
		\tilde  a_3  &= - \frac{C_\kappa^2  C_d^2 \sin^3\varphi \cos\varphi}{2}, \label{subeq-3: coef of aeff nd_x}\\
		\tilde  a_4  &= \frac{C_\kappa^2 C_d^2 \sin^4\varphi \cos^2\varphi}{16}. \label{subeq-4: coef of aeff nd_x}
	\end{align}

\end{subequations}

\section{CFD model for valve discharge}
\label{sec:Appendix_II}

The purpose of this section is to provide more details of the CFD can be used to calculate accurate discharge characteristics.
In truth, valve discharge involves complex flow regimes and
coefficients $C_{d}$ or characteristics $A_{\rm eff}(x)$ are just
approximations. To test the accuracy of such approximations, we used
the commercial code ANSYS-CFX to compute $A_\mathrm{eff}(x)$ and $C_d$
for some prototypical valve geometries.

\begin{figure}
\centering
\subfloat[]{
\includegraphics[width=0.45 \textwidth]{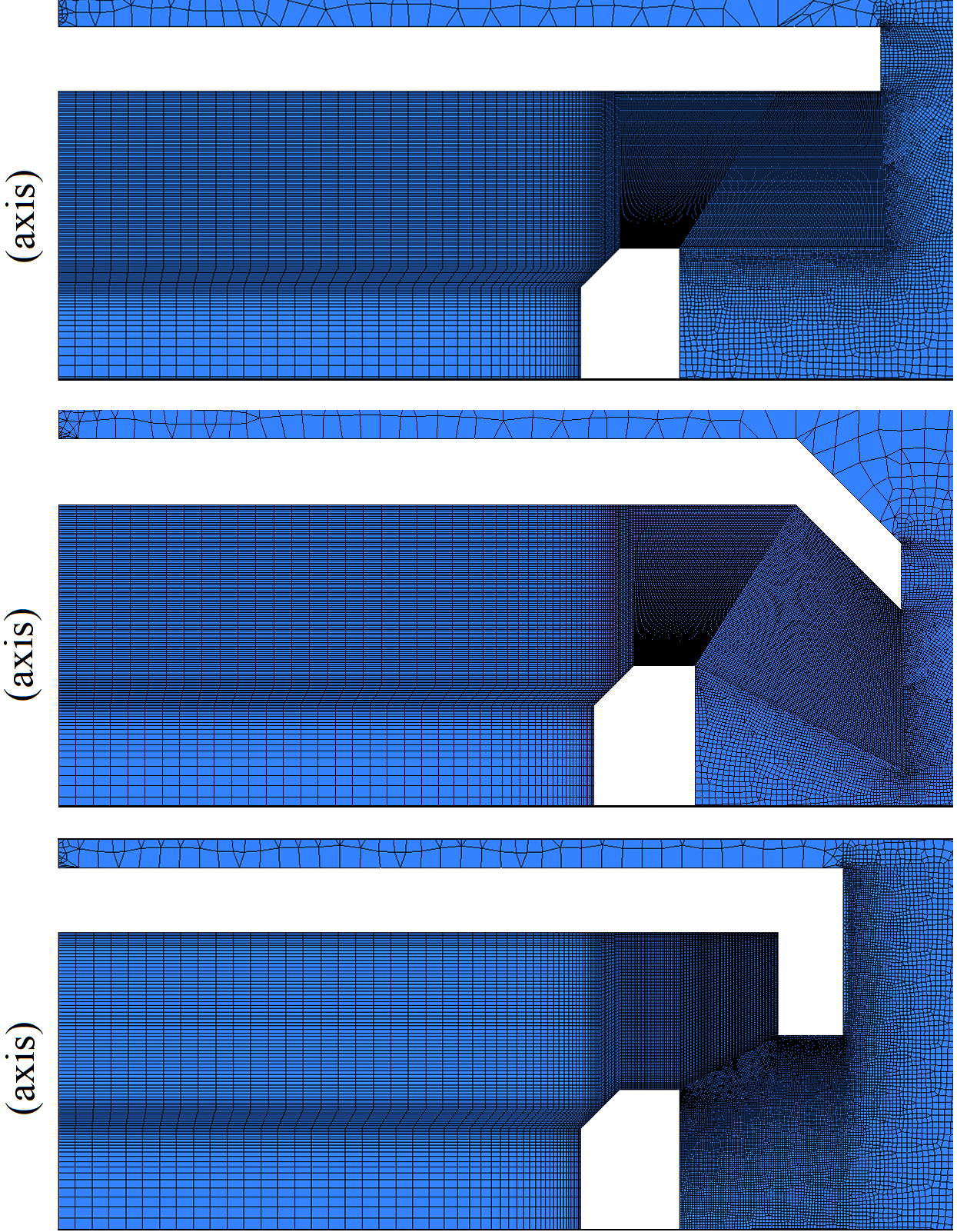}
\label{fig:Collar}
}
\subfloat[]{
\includegraphics[width=0.45 \textwidth]{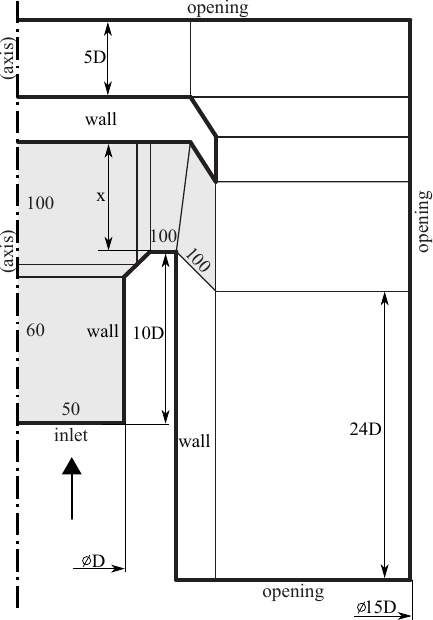}
\label{fig:CollarMeshBC}
}
\label{fig:AeffJetAngleGrid}
\caption{(a) 
The disc geometries of three disk-shaped valves (the $\varphi$ collar angles equal $90^\circ$, $45^\circ$ and $0^\circ$ from top to bottom) along with the numerical mesh. (b) The number of nodes and the boundary conditions for the middle case. The filled parts indicate the structured parts of the mesh, the unfilled regions include unstructured mesh. (The figure is
not proportional for illustrative purposes.)}
\end{figure}

The valve shapes and the corresponding meshes can be seen in Fig.~\ref{fig:Collar} (the mesh for the conical valve has a similar resolution,
see Fig.~\ref{figAeffJetAngleGrid} for its outline). 

Hybrid meshes consisting of both structured and unstructured parts
were created in GMSH for all geometries, with approximately 82\,000
elements for the disc-shaped bodies and 63\,000 for the conical
case. A $10D$ long inlet piping was added so that the inlet boundary
was sufficiently far from the valve. The total pressure at the pipe
inlet was varied between $1$ and $6$\,bar with $1 \,$ bar increments,
and the relative valve lift $x/D$ was varied between $0.2$ and $1.2$
($1.0$ for the cone) in increments of 0.2. In all cases, the ambient
pressure prescribed on the opening boundaries was $1$\,bar and the
inlet and ambient temperatures were $293\,$K. The boundary conditions
and the number of elements at the structured parts can be seen in
Fig. \ref{fig:CollarMeshBC}.

The calculations were performed with incompressible water and the $k-\epsilon$ turbulence model. The monitored quantities for post-processing were the mass flow rate, the force acting on the valve, and the static pressure upstream of the valve. 

The simulations were set up in ANSYS-CFX with the High Resolution
scheme for the advection terms and the turbulence model, and the
second-order backward Euler option for time stepping. A transient RANS
solver was used, from which the steady-state results were extracted
(the steady-state solver did not produce satisfactory convergence). A
mesh dependency study was conducted by doubling the mesh resolution on
all edges for two geometries (the $45^\circ$ disc and the conical cases) in
four seperate cases for each (largest and smallest lifts with highest and lowest inlet
pressures). We found that between the baseline and the high resolution
meshes, the relative differences in the fluid forces and the mass flow
rates were all below $1\%$ for the disc and
were below $3\%$ for the conical geometry 

In addition to the results in the main body, we report here results for
the computed discharge coefficients, which are summarized in Table \ref{tableCD} for the different geometries. We found that $C_d$ to be independent of the pressure difference within this parameter range but do depend somewhat on the valve lift.  In particular, 
the coefficient at higher lifts is significantly lower for the disc-shaped geometries. The reason for this observation could be that the flow-through area between the seat and the body is not as well-defined
as in the case of the conical body, for which we found the changes in the discharge coefficient to be much less significant.

\begin{table}
\centering
\begin{tabular}{| c || c | c | c | c |}
\hline
 & \multicolumn{3}{ c |}{Disc} & Cone \\
\hline
$x/x_{\rm max}$ & $C_{d,0^\circ}$ & $C_{d,45^\circ}$ & $C_{d,90^\circ}$ & $C_{d,135^\circ}$\\
\hline
$0.2$ & $0.9032$ & $0.9192$ & $0.9082$ & $0.9635$ \\
\hline
$0.4$ & $0.6989$ & $0.7835$ & $0.7655$ & $0.8722$ \\
\hline
$0.6$ & $0.6234$ & $0.7133$ & $0.7501$ & $0.8565$ \\
\hline
$0.8$ & $0.6070$ & $0.6785$ & $0.7070$ & $0.8520$ \\
\hline
$1.0$ & $0.5786$ & $0.6275$ & $0.6707$ & $0.8877$ \\
\hline
$1.2$ & $0.5496$ & $0.6064$ & $0.6553$ & --- \\
\hline
\end{tabular}
\caption{The discharge coefficient versus the relative lift for different discharge angles.}
\label{tableCD}
\end{table}


\section{Derivation of the quarter-wave model}\label{sec:Appendix_III}

The purpose of this section is to provide 
definitive details of the derivation of the  quarter-wave model in the presence of convective nonlinearity, pipe friction and inlet pressure loss. We take $N = 1$ in Eq.~\eqref{eq:trial-solution-coll-dimensional}, thus leading to pressure and velocity distributions given by
\vspace{-0.5em}
	\begin{align*}
		p(\xi,t) &=p_0(t) + B(t)\sin(\frac{\pi \xi}{2L}), \\
		v(\xi,t) &= v_L(t) + C(t)\cos(\frac{\pi \xi}{2L}), 
	\end{align*}
where $v_L$ is defined in Eq.~\eqref{eq:BC-V-at-1} as
\begin{equation}
v_L(t)= C_d C_\kappa \frac{A_{\rm ft} (x(t)}{A_0} \sqrt{\frac{p_0+B-p_b}{\rho}}  
\label{eq:vL_qw}
\end{equation}
in this context (without the term $p_b$ in the case of gas flow) 
and 
\begin{equation}
  p_0 = p_r \left [1- \chi(p_r,v_L+C) \right ].
  \label{eq:p0_qw}
  \end{equation}

We then substitute the assumed pressure and velocity distributions into the 
pipeline fluid
dynamics equations, and force them to be satisfied at the midpoint of the pipe, i.e.~at $\xi = L/2$.
Then we can arrive at
\begin{align}
	\sqrt { 2 } \dot { p } _ 0 + \dot { B } + ( \sqrt { 2 } v _ { L } + C ) \frac { \pi } { 2\sqrt{2} L } B - a ^ {
		2 } \rho  \frac { \pi } { 2 L } C &= 0 \\
	\sqrt { 2 } \dot { v } _ { L } + \dot { C } - ( \sqrt { 2 } v _ { L } + C ) \frac { \pi } { 2 \sqrt{2}L } C +
	\frac { 1 } { \rho } \frac { \pi } { 2 L } B & =  {-}\lambda \frac { L } { D } ( \sqrt { 2 } v_L+ C )
	\left |\sqrt { 2 } v_L + C \right|
\end{align}

In the end, the whole system is approximated by five ODEs written in the implicit form
\begin{subequations}
	\label{eq: OEQ of QWM}
	\begin{align}
		m\ddot{x} & = (p_0 + B - p_b) {A}_{\rm eff}(x) - c \dot{x}  - k(x + x_{\rm pre})\\ 
		\dot{p}_r & = \frac{a^2}{V} [ \dot{m}_{\rm in} - A_0\rho (v_L(t) + C(t)) ]\\
		\dot{B} &= - \sqrt{2} \dot{p}_0 + a ^ { 2 } \rho  \frac { \pi } { 2 L }C - (\sqrt{2} v_{L} + C) \frac
		{\pi}{2\sqrt{2}L} B\\
		\dot{C} &=  -\sqrt {2}\dot{v}_L + ( \sqrt { 2 } v_{ L } + C ) \frac { \pi } { 2 \sqrt{2}L } C -\frac { 1 } {
			\rho } \frac { \pi } { 2 L } B  {-}\lambda \frac { L } { D } ( \sqrt { 2 } v_L+ C ) \left | \sqrt {
			2 } v_L + C \right |
	\end{align}
\end{subequations}
for the dynamic variables $\mathbf{x} = (x,\dot{x},p_r,B,C)^T$.
The equations can be closed by differentiation of the relations Eqs.~\eqref{eq:p0_qw} and \eqref{eq:vL_qw} to write 
$$
\dot{p}_0  = \dot{p}_r -\frac{\partial \chi }{ \partial p_r } \dot{p}_r - \frac{\partial \chi }{ \partial v} (\dot{v}_L + C)
$$
and 
$$
\dot{v}_L = \frac{C_d C_\kappa}{\sqrt{\rho} A_0} 
\left ( A^\prime(x)_{\rm ft} \sqrt{p_0+B-p_b} \dot{x} 
- A_{\rm ft}(x) \frac{\dot{p}_0 + \dot{B}}{\sqrt{p_0+B-p_b}} 
\right )
$$

Hence the system Eq.~\eqref{eq: OEQ of QWM} can be written in the form 
$$
M(\mathbf{x}) \dot{\mathbf{x}} = \mathbf{F}(\mathbf{x})
$$
where $M$ is a $5\times 5$ matrix which is non-singular for 
sufficiently small $\mathbf{x}$ and $\mathbf{F}$ is a 5-dimensional vector of its arguments.

Given an equilibrium of the form 
$$
\mathbf{x}_e = \begin{pmatrix}
x_e \\
0 \\
P(x_e) \\
0 \\
0
\end{pmatrix}
$$
then the Jacobian matrix of the equations at this point
can be written in the form 
$$
J(\mathbf{x}_e) = \left .\frac{ \partial }{\partial \mathbf{x}} M^{-1} \mathbf{F}
\right |_{\mathbf{x}=\mathbf{x}_e} .
$$
It is the eigenvalues of this matrix that determine the stability of the quarter-wave model. The corresponding expressions are quite convoluted and are best investigated using computer algebra. Detailed expressions are available in \cite{tang_thesis}.